\documentclass[a4paper]{article}

\usepackage{amsmath}
\usepackage{amssymb}
\usepackage{amsthm}
\usepackage{mathtools}
\usepackage{tikz}
\usepackage{tikz-cd}
\usepackage{geometry}
\usepackage{hyperref}
\usepackage{booktabs}
\usepackage{bussproofs}
\usepackage{cleveref}
\usepackage[shortlabels]{enumitem}
\usepackage{subcaption}
\usepackage{xspace}

\usetikzlibrary{fit}

\newcommand\ord[1]{\mathbf{#1}}

\newcommand\cat[1]{\mathbf{#1}}

\newcommand\Cat{\mathbf{Cat}}
\newcommand\Pos{\mathbf{Pos}}
\newcommand\Set{\mathbf{Set}}
\newcommand\Ord{\mathbf{\Delta}}
\newcommand\Int{\boldsymbol{\nabla}}

\newcommand\email[1]{\href{mailto:#1}{#1}}

\newcommand\hio{\textsc{homotopy.io}\xspace}

\DeclarePairedDelimiter\set{\{}{\}}
\DeclarePairedDelimiter\abs{\lvert}{\rvert}

\DeclareMathOperator\id{id}
\DeclareMathOperator\colim{colim}

\DeclareMathOperator\Pt{\mathsf{Pt}}

\DeclareMathOperator\Zig{\mathsf{Zig}}
\DeclareMathOperator\Exp{\mathsf{Exp}}

\DeclareMathOperator\Reg{\mathsf{R}}
\DeclareMathOperator\reg{\mathsf{r}}
\DeclareMathOperator\sing{\mathsf{s}}

\newtheorem{theorem}{Theorem}

\newtheorem{proposition}[theorem]{Proposition}

\theoremstyle{definition}
\newtheorem{definition}[theorem]{Definition}
\newtheorem{example}[theorem]{Example}

\title{A layout algorithm for higher-dimensional string diagrams}
\author{%
Calin Tataru
\thanks{University of Cambridge, \email{calin.tataru@cl.cam.ac.uk}}
\and
Jamie Vicary
\thanks{University of Cambridge, \email{jamie.vicary@cl.cam.ac.uk}}
}

\begin{document}

\maketitle

\begin{abstract}
The algebraic zigzag construction has recently been introduced as a combinatorial foundation for a higher dimensional notion of string diagram.
For use in a proof assistant, a layout algorithm is required to determine the optimal rendering coordinates, across multiple projection schemes including 2D, 3D, and 4D.
For construction of these layouts, a key requirement is to determine the linear constraints which the geometrical elements must satisfy in each dimension.
Here we introduce a new categorical tool called \emph{injectification}, which lifts a functorial factorisation system on a category to diagrams over that category, and we show that implementing this recursively in the category of finite posets allows us to systematically generate the necessary constraints.
These ideas have been implemented as the layout engine of the proof assistant \hio, enabling attractive and practical visualisations of complex higher categorical objects.
\end{abstract}

\section{Introduction} \label{sec:introduction}

Associative $n$-categories have been recently proposed as a new semistrict model of $n$-categories with strictly associative and unital composition in all dimensions~\cite{Dor18}.
They retain sufficient weak structure, in the form of homotopies of composites, to reasonably conjecture that they are equivalent to weak $n$-categories.
This makes them an attractive language for a proof assistant, because they allow composites to be unique without sacrificing expressiveness.
Indeed, the theory is the basis for the proof assistant \hio which allows users to construct and manipulate terms in finitely-presented $n$-categories.

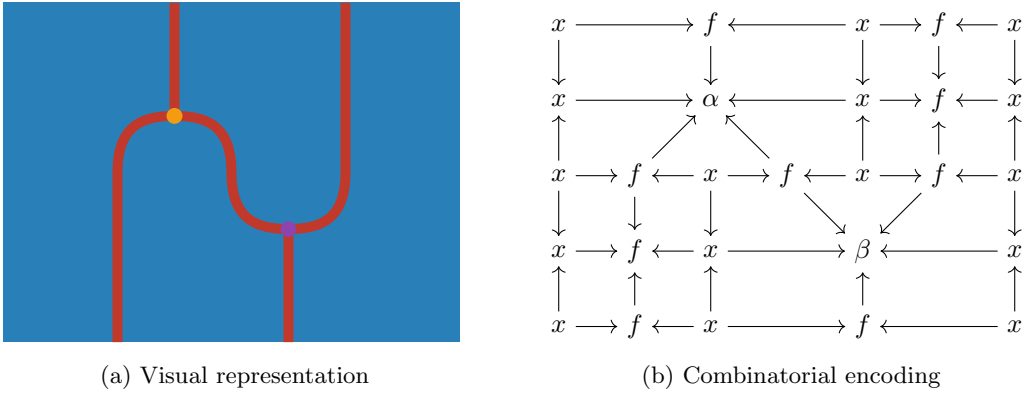
\begin{figure}
\centering
\begin{subfigure}{0.49\textwidth}
\centering
\scalebox{0.75}{\begin{tikzpicture}
\definecolor{generator-3-2-0-pos}{RGB}{142, 68, 173}
\definecolor{generator-1-1-0-pos}{RGB}{192, 57, 43}
\definecolor{generator-2-2-0-pos}{RGB}{243, 156, 18}
\definecolor{generator-0-0-0-pos}{RGB}{41, 128, 185}
\begin{scope}
% Background surfaces
\fill[generator-0-0-0-pos] (0,0) -- (8,0) -- (8,6) -- (0,6) -- (0,0);
% Wire layers
\draw[color=generator-1-1-0-pos, line width=5pt](2,0) -- (2,3) .. controls (2,3.8) and (2.4,4) .. (3,4) -- (3,6)(5,0) -- (5,2) .. controls (4.4,2) and (4,2.2) .. (4,3) .. controls (4,3.8) and (3.6,4) .. (3,4)(5,2) .. controls (5.6,2) and (6,2.2) .. (6,3) -- (6,6);
\end{scope}
\fill[generator-3-2-0-pos] (5,2) circle (0.14);
\fill[generator-2-2-0-pos] (3,4) circle (0.14);
\end{tikzpicture}}
\caption{Visual representation}
\end{subfigure}
\begin{subfigure}{0.49\textwidth}
\centering
\begin{tikzpicture}

\node (R0R0) at (0, 0) {$x$};
\node (R0S0) at (1, 0) {$f$};
\node (R0R1) at (2, 0) {$x$};
\node (R0S1) at (4, 0) {$f$};
\node (R0R2) at (6, 0) {$x$};

\node (S0R0) at (0, 1) {$x$};
\node (S0S0) at (1, 1) {$f$};
\node (S0R1) at (2, 1) {$x$};
\node (S0S1) at (4, 1) {$\beta$};
\node (S0R2) at (6, 1) {$x$};   

\node (R1R0) at (0, 2) {$x$};
\node (R1S0) at (1, 2) {$f$};
\node (R1R1) at (2, 2) {$x$};
\node (R1S1) at (3, 2) {$f$};
\node (R1R2) at (4, 2) {$x$};
\node (R1S2) at (5, 2) {$f$};
\node (R1R3) at (6, 2) {$x$};

\node (S1R0) at (0, 3) {$x$};
\node (S1S0) at (2, 3) {$\alpha$};
\node (S1R1) at (4, 3) {$x$};
\node (S1S1) at (5, 3) {$f$};
\node (S1R2) at (6, 3) {$x$};

\node (R2R0) at (0, 4) {$x$};
\node (R2S0) at (2, 4) {$f$};
\node (R2R1) at (4, 4) {$x$};
\node (R2S1) at (5, 4) {$f$};
\node (R2R2) at (6, 4) {$x$};

\draw[->] (R0R0) to (R0S0);
\draw[->] (R0R1) to (R0S0);
\draw[->] (R0R1) to (R0S1);
\draw[->] (R0R2) to (R0S1);

\draw[->] (S0R0) to (S0S0);
\draw[->] (S0R1) to (S0S0);
\draw[->] (S0R1) to (S0S1);
\draw[->] (S0R2) to (S0S1);

\draw[->] (R1R0) to (R1S0);
\draw[->] (R1R1) to (R1S0);
\draw[->] (R1R1) to (R1S1);
\draw[->] (R1R2) to (R1S1);
\draw[->] (R1R2) to (R1S2);
\draw[->] (R1R3) to (R1S2);

\draw[->] (S1R0) to (S1S0);
\draw[->] (S1R1) to (S1S0);
\draw[->] (S1R1) to (S1S1);
\draw[->] (S1R2) to (S1S1);

\draw[->] (R2R0) to (R2S0);
\draw[->] (R2R1) to (R2S0);
\draw[->] (R2R1) to (R2S1);
\draw[->] (R2R2) to (R2S1);

\draw[->] (R0R0) to (S0R0);
\draw[->] (R0R1) to (S0R1);
\draw[->] (R0R2) to (S0R2);

\draw[->] (R0S0) to (S0S0);
\draw[->] (R0S1) to (S0S1);

\draw[->] (R1R0) to (S0R0);
\draw[->] (R1R1) to (S0R1);
\draw[->] (R1R3) to (S0R2);

\draw[->] (R1S0) to (S0S0);
\draw[->] (R1S1) to (S0S1);
\draw[->] (R1S2) to (S0S1);

\draw[->] (R1R0) to (S1R0);
\draw[->] (R1R2) to (S1R1);
\draw[->] (R1R3) to (S1R2);

\draw[->] (R1S0) to (S1S0);
\draw[->] (R1S1) to (S1S0);
\draw[->] (R1S2) to (S1S1);

\draw[->] (R2R0) to (S1R0);
\draw[->] (R2R1) to (S1R1);
\draw[->] (R2R2) to (S1R2);

\draw[->] (R2S0) to (S1S0);
\draw[->] (R2S1) to (S1S1);

\end{tikzpicture}
\caption{Combinatorial encoding}
\end{subfigure}
\caption{Example of a 2-diagram in \hio.}
\label{fig:example}
\end{figure}

Terms in the theory, called \emph{$n$-diagrams}, are inherently geometrical and have a direct interpretation as $n$-dimensional string diagrams.
The proof assistant \hio embraces this, and the only way to view or interact with terms is via a visual representation.
However, this poses an obstacle, because $n$-diagrams are combinatorial objects, defined via an inductive construction, which do not contain \emph{layout} information indicating how they should be drawn.

\begin{figure}
\centering
\includegraphics[width=0.5\textwidth]{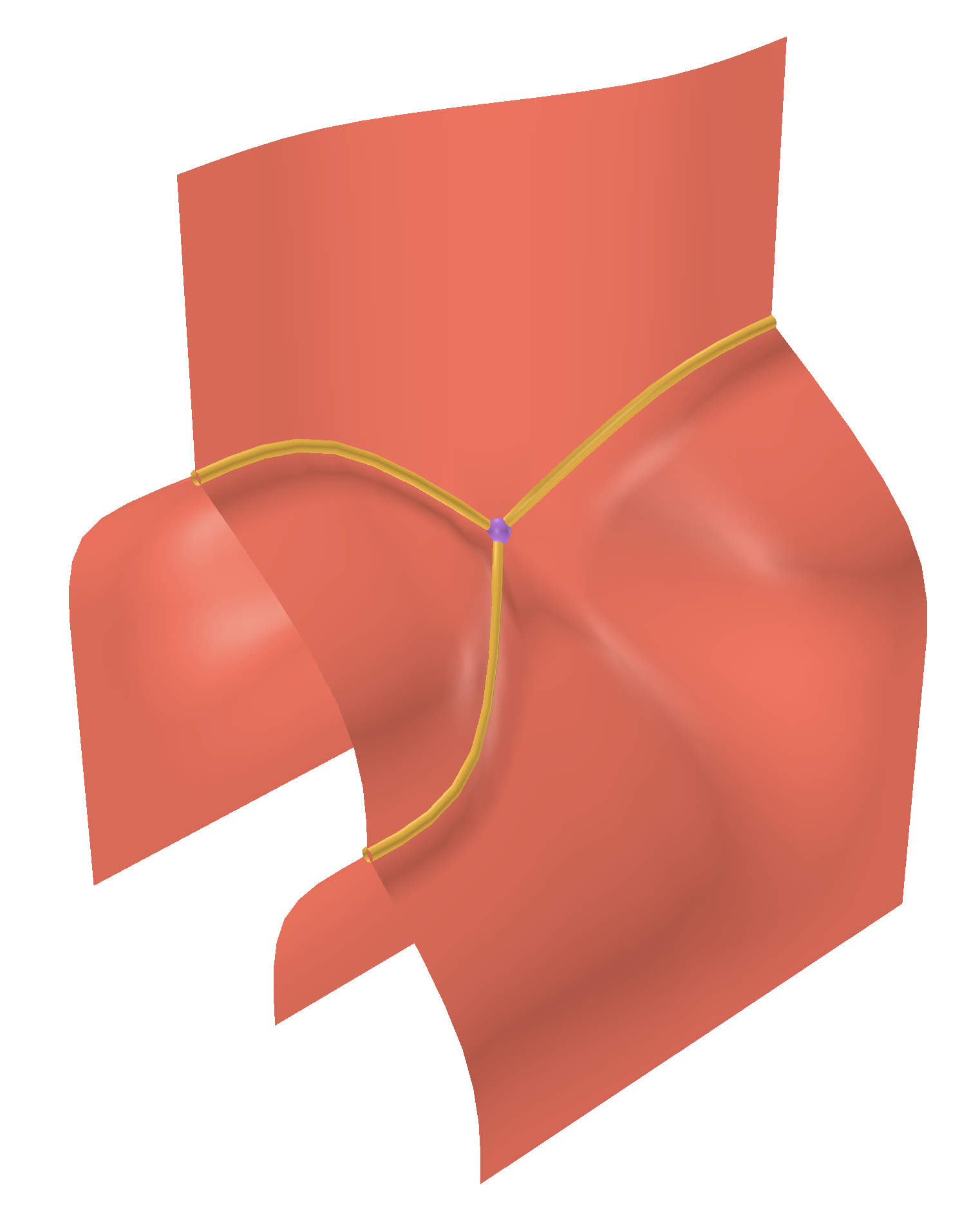}
\caption{Example of a 3-diagram in \hio.}
\label{fig:associator}
\end{figure}

Constructing a visually appealing layout is a non-trivial problem, particularly in higher dimensions.
Previous proof assistants have used ad-hoc methods~\cite{BKV18}, but these are inefficient and become infeasible in higher dimensions, so a more general solution is needed.
Our approach is to use a functorial factorisation system to produce a set of constraints, which are then passed to a linear solver, which constructs the final layout.
The algorithm has been successfully implemented in the new version of \hio\footnotemark{}, where it can construct the following output:
\begin{itemize}
\item In dimension~2, it yields 2D string diagrams (see \Cref{fig:example}).
\item In dimension~3, it yields 3D surface diagrams (see \Cref{fig:associator}).
\item In dimension~4, it yields smooth animations of 3D surface diagrams (see \cref{sec:examples}).
\end{itemize}
In all these cases, the algorithm produces high-quality output which improves the ability of users to comprehend complex diagrams in higher dimensions.
While ad-hoc approaches would likely suffice in dimension 2, and with difficulty in dimension 3, the abstract approach we develop here becomes essential from dimension 4 onwards, with the constraints becoming highly complex.

\footnotetext{Currently in beta stage, the new version of the tool is hosted at \url{beta.homotopy.io}.} 

Our scheme is motivated as follows.
The $n$-diagrams in \hio are defined over the category $\Ord$ of finite total orders and monotone functions, using the zigzag construction, which we review in \cref{sec:zigzags}; this construction technique is inductive on the dimension of the diagram.
In some sense, this total order data encodes all linear relationships between parts of a diagram: for example, in a 2D diagram, that one wire is to the left of another; or in a 3D diagram, that one sheet is in front of another.

To obtain the full set of layout constraints, a robust intuition is that we should take the \emph{colimit} of all this total order data, yielding a single combinatorial object that contains all the constraints.
However, this approach does not immediately succeed, for a number of reasons.
Firstly, the category $\Ord$ is logically poor, lacking the necessary colimit objects.
Secondly, the non-injectivity of some morphisms of $\Ord$ can cause parts of a diagram to collapse to a point; for an attractive layout this should be prevented, otherwise we would have a degenerate layout.
Thirdly, the inductive nature of our construction means that we need to gather constraints recursively, across all dimensions of the given $n$-diagram. 

Here we provide new logical tools which solve these problems, giving a new high-level layout algorithm which works in all dimensions.
Regarding the problems above, we overcome the first by passing from $\Ord$ to $\Pos$, the category of finite posets and monotone functions; this is technically straightforward since $\Ord$ is a full subcategory $\Pos$, and since $\Pos$ has all finite colimits, our construction can proceed.
The issue of non-injectivity is addressed by exploiting a mono-epi factorisation system in $\Pos$, which allows us to realize our diagram structures entirely via monomorphisms, a process we call \emph{injectification}.
Finally, we handle the inductive structure by showing that this injectification process is functorial, allowing the different dimensions of our diagram to be handled uniformly.

To illustrate how these problems are solved, consider the following diagram in $\Ord$:
\[
\begin{tikzpicture}
\node (a) at (1.5, 0) {$a$};
\node (b) at (2.5, 0) {$b$};
\node (c) at (0, 1) {$c$};
\node (d) at (1, 1) {$d$};
\node (e) at (3, 1) {$e$};
\node (f) at (4, 1) {$f$};

\node[draw,dotted,rounded corners,fit=(a) (b)] {};
\node[draw,dotted,rounded corners,fit=(c) (d)] {};
\node[draw,dotted,rounded corners,fit=(e) (f)] {};

\draw[->] (a) to (b);
\draw[->] (c) to (d);
\draw[->] (e) to (f);

\draw[->] (a) to (c);
\draw[->] (b) to (c);

\draw[->] (a) to (e);
\draw[->] (b) to (e);
\end{tikzpicture}
\]
This does \emph{not} have a colimit in $\Ord$, but this is solved by passing to $\Pos$, where the colimit exists and is given by:
\[
\begin{tikzpicture}
\node (a) at (1.5, 0) {$a$};
\node (b) at (2.5, 0) {$b$};
\node (c) at (0, 1) {$c$};
\node (d) at (1, 1) {$d$};
\node (e) at (3, 1) {$e$};
\node (f) at (4, 1) {$f$};

\node[draw,dotted,rounded corners,fit=(a) (b)] {};
\node[draw,dotted,rounded corners,fit=(c) (d)] {};
\node[draw,dotted,rounded corners,fit=(e) (f)] {};

\draw[->] (a) to (b);
\draw[->] (c) to (d);
\draw[->] (e) to (f);

\draw[->] (a) to (c);
\draw[->] (b) to (c);

\draw[->] (a) to (e);
\draw[->] (b) to (e);

\node (g) at (1.5, 3) {$g$};
\node (h) at (2.5, 2.5) {$h$};
\node (i) at (2.5, 3.5) {$i$};

\node[draw,dotted,rounded corners,fit=(g) (h) (i)] {};

\draw[->] (g) to (h);
\draw[->] (g) to (i);

\draw[->] (c) to (g);
\draw[->] (d) to (h);
\draw[->] (e) to (g);
\draw[->] (f) to (i);
\end{tikzpicture}
\]
However, the colimit forces $a$ and $b$ to be identified, which is undesirable for a layout algorithm.
Injectification solves this by replacing the diagram with one where every map is injective, such as:
\[
\begin{tikzpicture}
\node (a) at (2.5, 0) {$a$};
\node (b) at (3.5, 0) {$b$};
\node (c1) at (0, 1) {$c_1$};
\node (c2) at (1, 1) {$c_2$};
\node (d) at (2, 1) {$d$};
\node (e1) at (4, 1) {$e_1$};
\node (e2) at (5, 1) {$e_2$};
\node (f) at (6, 1) {$f$};

\node[draw,dotted,rounded corners,fit=(a) (b)] {};
\node[draw,dotted,rounded corners,fit=(c1) (c2) (d)] {};
\node[draw,dotted,rounded corners,fit=(e1) (e2) (f)] {};

\draw[->] (a) to (b);
\draw[->] (c1) to (c2);
\draw[->] (c2) to (d);
\draw[->] (e1) to (e2);
\draw[->] (e2) to (f);

\draw[->] (a) to (c1);
\draw[->] (b) to (c2);

\draw[->] (a) to (e1);
\draw[->] (b) to (e2);

\node (g1) at (2, 3) {$g_1$};
\node (g2) at (3, 3) {$g_2$};
\node (h) at (4, 2.5) {$h$};
\node (i) at (4, 3.5) {$i$};

\node[draw,dotted,rounded corners,fit=(g1) (g2) (h) (i)] {};

\draw[->] (g1) to (g2);
\draw[->] (g2) to (h);
\draw[->] (g2) to (i);

\draw[->] (c1) to (g1);
\draw[->] (c2) to (g2);
\draw[->] (d) to (h);
\draw[->] (e1) to (g1);
\draw[->] (e2) to (g2);
\draw[->] (f) to (i);
\end{tikzpicture}
\]
This colimit encodes the correct set of constraints, that can then be passed to a linear solver to generate the actual layout data.
In this case the problem is of course trivial to solve by hand, but in realistic diagrams produced in the proof assistant, particularly in high dimensions, the resulting posets can have thousands of elements.

We prove correctness of our algorithm in an abstract way, entirely via categorical methods, and our results could therefore be generalized to other sorts of layout problems by changing the underlying categories and factorisation systems that are employed.

\subsection{Related work}

The theory of associative $n$-categories was originally developed by Dorn, Douglas, and Vicary, and was described in Dorn's thesis~\cite{Dor18}.
The zigzag construction, which gives an inductive definition of terms in associative $n$-categories, was described by Reutter and Vicary~\cite{RV19} for the purpose of describing a contraction algorithm for building non-trivial terms in the theory.
This was followed by Heidemann, Reutter, and Vicary's paper on a normalisation algorithm~\cite{HRV22}, which allows type checking in the proof assistant.

There are other proof assistants that use string diagrams as a visualisation method, and these all use a layout algorithm to display their results (see \Cref{tab:layout}).
Of these, \textsc{wiggle.py}~\cite{Bur23} is perhaps the most similar, since that system also uses linear programming compute the layout; that system is a dedicated 3D renderer and generates the constraints via custom code, whereas we generate the constraints in any number of dimensions via a uniform mathematical approach.

\begin{table}
\centering
\begin{tabular}{@{}llll@{}}
\toprule
proof assistant                          & generality                    & dimensions & layout algorithm          \\ \midrule
\textsc{Cartographer}~\cite{SWZ19}       & symmetric monoidal categories & two        & layered graph drawing     \\
\textsc{Quantomatic}~\cite{KZ15}         & compact closed categories     & two        & force-based graph drawing \\
\textsc{wiggle.py}~\cite{Bur23}          & monoidal bicategories         & three      & linear programming        \\ \bottomrule \\
\end{tabular}
\caption{Comparison of layout algorithms in other proof assistants.}
\label{tab:layout}
\end{table}

\subsection{Notation}

For each $n \in \mathbb{N}$, we write $\ord n$ for the totally-ordered set $\set{0, \dots, n - 1}$.
We write $\Ord$ for the category where objects are natural numbers $n \in \mathbb{N}$ and morphisms are monotone maps $\ord n \to \ord m$.
We also write $\Int$ for the subcategory of $\Ord$ with non-zero natural numbers and monotone maps preserving the first and last elements.
Finally, we write $\Set$ for the category of sets, and $\Pos$ for the category of finite partial orders (or posets).

\subsection{Acknowledgements}

The authors are grateful to Nick Hu, Ioannis Markakis, and Chiara Sarti for useful discussions, and Nathan Corbyn for implementing the WebGL renderer in \hio which allows the layout algorithm to be used in practice.

\section{Zigzag construction} \label{sec:zigzags}

We begin by recalling the zigzag construction developed by Reutter and Vicary~\cite{RV19}.
This gives a simple inductive definition of \emph{$n$-diagrams}, which form the terms of the theory of associative $n$-categories.

\begin{definition}
In a category $\cat C$, a \emph{zigzag} $X$ is a finite sequence of cospans:
\[
\begin{tikzcd}[column sep=small]
X(\reg 0) \rar & X(\sing 0) & X(\reg 1) \lar \rar & X(\sing 1) & \cdots \lar \rar & X(\sing n - 1) & X(\reg n) \lar
\end{tikzcd}
\]
This consists of \emph{singular objects} of the form $X(\sing i)$ for $i \in \ord n$, and \emph{regular objects} of the form $X(\reg j)$ for $j \in \ord{n + 1}$.
Such a zigzag is said to have \emph{length} $n$, given by the number of singular objects.
We write $X_\mathsf{s} = \ord n$ and $X_\mathsf{r} = \ord{n + 1}$.
\end{definition}

There is a natural notion of morphism between zigzags.
It is based on an equivalence of categories 
\[
\Reg : \Ord \xrightarrow{\sim} \Int^\mathrm{op}
\]
first observed by Wraith~\cite{Wra93}, sending $n$ to $n + 1$ and each monotone map $f : \ord n \to \ord m$ in $\Ord$ to the monotone map $\Reg f : \ord{m + 1} \to \ord{n + 1}$ in $\Int$, given by $(\Reg f)(i) \coloneqq \min(\set{j \in \ord n \mid f(j) \geq i} \cup \set{n})$.
The equivalence has a nice geometrical intuition, which is illustrated in \Cref{fig:interleaving}: the map $f$ going up the page is ``interleaved'' with the map $\Reg f$ going down the page.
We can obtain the definition of a zigzag map from this picture by labelling the arrows with actual morphisms in $\cat C$ (for the arrows of $f$) or $\cat C^\mathrm{op}$ (for the arrows of $\Reg f$).

\begin{figure}
\centering
\begin{tikzcd}[column sep=small]
\color{blue}{\bullet} \ar[dr, blue] & \color{red}{\bullet} & \color{blue}{\bullet} \ar[dr, blue] & \color{red}{\bullet} & \color{blue}{\bullet} \ar[dl, blue] & \color{red}{\bullet} & \color{blue}{\bullet} \ar[dr, blue] & \color{red}{\bullet} & \color{blue}{\bullet} \ar[dl, blue] \\
& \color{blue}{\bullet} & \color{red}{\bullet} \ar[ul, red] & \color{blue}{\bullet} & \color{red}{\bullet} \ar[ur, red] & \color{blue}{\bullet} & \color{red}{\bullet} \ar[ul, red] & \color{blue}{\bullet}
\end{tikzcd}
\caption{An interleaving of $\color{red}f : \ord 3 \to \ord 4$ in $\Ord$ going up, and $\color{blue}\Reg f : \ord 5 \to \ord 4$ in $\Int$ going down.}
\label{fig:interleaving}
\end{figure}

\begin{definition}
A \emph{zigzag map} $f : X \to Y$ between a zigzag $X$ of length $n$ and a zigzag $Y$ of length $m$ consists of the following data:
\begin{itemize}
\item a \emph{singular map} in $\Ord$,
\[
f_\mathsf{s} : \ord n \to \ord m
\]
with an implied \emph{regular map} in $\Int$,
\[
f_\mathsf{r} \coloneqq \Reg f_\mathsf{s} : \ord{m + 1} \to \ord{n + 1}
\]
\item for every $i \in \ord n$, a \emph{singular slice} in $\cat C$,
\[
f(\sing i) : X(\sing i) \to Y(\sing f_\mathsf{s}(i))
\]
\item for every $j \in \ord{m + 1}$, a \emph{regular slice} in $\cat C$,
\[
f(\reg j) : X(\reg f_\mathsf{r}(j)) \to Y(\reg j)
\]
\end{itemize}
These must satisfy the following conditions, for every $i \in \ord m$:
\begin{itemize}
\item if $f_\mathsf{s}^{-1}(i)$ is empty, this diagram must commute:
\[
\begin{tikzcd}[column sep=small]
Y(\reg i) \rar & Y(\sing i) & Y(\reg i + 1) \lar \\
& X(\reg f_\mathsf{r}(i)) \ular["f(\reg i)"] \urar["f(\reg i + 1)"'] &
\end{tikzcd}
\]
\item if $f_\mathsf{s}^{-1}(i)$ is non-empty, the following diagrams must commute, where $p$ and $q$ are the lower and upper bounds of the preimage, and $p \leq j < q$:
\[
\begin{tikzcd}[column sep=small]
Y(\reg i) \rar & Y(\sing i) \\
X(\reg p) \rar \uar["f(\reg i)"] & X(\sing p) \uar["f(\sing p)"']
\end{tikzcd}
\quad
\begin{tikzcd}[column sep=small]
& Y(\sing i) & \\
X(\sing j) \urar["f(\sing j)"] & X(\reg j + 1) \lar \rar & X(\sing j + 1) \ular["f(\sing j + 1)"']
\end{tikzcd}
\quad
\begin{tikzcd}[column sep=small]
Y(\sing i) & Y(\reg i + 1) \lar \\
X(\sing q) \uar["f(\sing q)"] & X(\reg q + 1) \lar \uar["f(\reg i + 1)"']
\end{tikzcd}
\]
\end{itemize}
We write $f(\reg j, \sing i) : X(\reg j) \to Y(\sing j)$ for the unique diagonal fillers in each of these cases, for $i \in \ord m$ and $j \in f_\mathsf{r}([i, i + 1])$.
We call these the \emph{diagonal slices}.
\end{definition}

\begin{figure}
\centering
\begin{tikzcd}[cramped, column sep=small]
Y(\reg 0) \rar & Y(\sing 0) & Y(\reg 1) \lar \rar & Y(\sing 1) & Y(\reg 2) \lar \rar & Y(\sing 2) & Y(\reg 3) \lar \rar & Y(\sing 3) & Y(\reg 4) \lar \\
& X(\reg 0) \rar \ar[ul, blue] & X(\sing 0) \ar[ul, red] & X(\reg 1) \lar \rar \ar[ul, blue] \ar[ur, blue] & X(\sing 1) \ar[ur, red] & X(\reg 2) \lar \rar & X(\sing 2) \ar[ul, red] & X(\reg 3) \lar \ar[ul, blue] \ar[ur, blue]
\end{tikzcd}
\caption{A zigzag map $f : X \to Y$ with the singular slices in red and regular slices in blue.}
\end{figure}
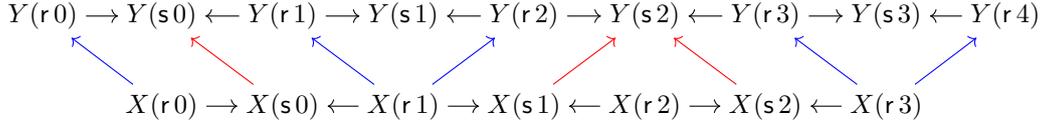

Given a category $\cat C$, we define its \emph{zigzag category} $\Zig(\cat C)$ to be the category whose objects are zigzags and whose morphisms are zigzag maps in $\cat C$.
The construction is functorial in $\cat C$, so we have a functor
\[
\Zig : \Cat \to \Cat
\]
The construction can be iterated, and we denote the \emph{$n$-fold zigzag category} by $\Zig^n(\cat C)$.

Note that $\Zig(\ord 1)$ is isomorphic to $\Ord$, since a zigzag in the terminal category is uniquely determined by its length.
Hence, since every category $\cat C$ has a unique functor $\cat C \to \ord 1$, this naturally induces a functor
\[
(-)_\mathsf{s} : \Zig(\cat C) \to \Ord
\]
called the \emph{singular projection}, which sends a zigzag to its length, and a zigzag map to its singular map.

\subsection{Diagrams}

Let $\Sigma$ be a set of labels equipped with a \emph{dimension} function $\dim : \Sigma \to \mathbb{N}$.
This naturally induces a partial order on $\Sigma$ given by $f < g$ iff $\dim f < \dim g$.
Therefore, we can view $\Sigma$ as a thin category and consider iterated zigzags over it.

\begin{definition}
A \emph{$\Sigma$-typed $n$-diagram} is an object of the $n$-fold zigzag category $\Zig^n(\Sigma)$.
\end{definition}

We think of the objects of $\Sigma$ as encoding generators in a signature.
Then the objects of $\Zig^n(\Sigma)$ can be seen as combinatorial encodings of $n$-dimensional string diagrams in the free $n$-category generated by that signature.
Not all $n$-diagrams are valid, however, and the full theory of associative $n$-categories has a \emph{typechecking} procedure for checking whether an $n$-diagram is \emph{well-typed} with respect to a signature (see \cite[Chapter 8]{Dor18} and \cite[Section 7]{HRV22}).

\begin{example}
\label{ex:monoid}
In the associative 2-category generated by a 0-cell $x$, a 1-cell $f : x \to x$, and a 2-cell $m : f \circ f \to f$, the 2-cell $m$ can be interpreted as a monad structure, and is encoded by the following 2-diagram:
\[
\begin{tikzpicture}

\node (R0R0) at (0, 0) {$x$};
\node (R0S0) at (1, 0) {$f$};
\node (R0R1) at (2, 0) {$x$};
\node (R0S1) at (3, 0) {$f$};
\node (R0R2) at (4, 0) {$x$};

\node (S0R0) at (0, 1) {$x$};
\node (S0S0) at (2, 1) {$m$};
\node (S0R1) at (4, 1) {$x$};

\node (R1R0) at (0, 2) {$x$};
\node (R1S0) at (2, 2) {$f$};
\node (R1R1) at (4, 2) {$x$};

\draw[->] (R0R0) to (R0S0);
\draw[->] (R0R1) to (R0S0);
\draw[->] (R0R1) to (R0S1);
\draw[->] (R0R2) to (R0S1);

\draw[->] (S0R0) to (S0S0);
\draw[->] (S0R1) to (S0S0);

\draw[->] (R1R0) to (R1S0);
\draw[->] (R1R1) to (R1S0);

\draw[->] (R0R0) to (S0R0);
\draw[->] (R0R2) to (S0R1);

\draw[->] (R0S0) to (S0S0);
\draw[->] (R0S1) to (S0S0);

\draw[->] (R1R0) to (S0R0);
\draw[->] (R1R1) to (S0R1);

\draw[->] (R1S0) to (S0S0);

\end{tikzpicture}
\]
This is defined over the signature $\Sigma = \set{x, f, m}$ with $\dim x = 0$, $\dim f = 1$, and $\dim m = 2$.
\end{example}

\begin{example}
In the associative 3-category generated by a 0-cell $x$ and a pair of 2-cells $\alpha, \beta : \id_x \to \id_x$, the \emph{Eckmann-Hilton move} is a 3-cell that can be interpreted as ``braiding $\alpha$ over $\beta$'', and is encoded by the following 3-diagram:
\[
\begin{tikzpicture}

% Source

\node (R0R0R0) at (1, 0) {$x$};

\node (R0S0R0) at (0, 1) {$x$};
\node (R0S0S0) at (1, 1) {$\alpha$};
\node (R0S0R1) at (2, 1) {$x$};

\node (R0R1R0) at (1, 2) {$x$};

\node (R0S1R0) at (0, 3) {$x$};
\node (R0S1S0) at (1, 3) {$\beta$};
\node (R0S1R1) at (2, 3) {$x$};

\node (R0R2R0) at (1, 4) {$x$};

\draw[->] (R0R0R0) to (R0S0R0);
\draw[->] (R0R0R0) to (R0S0R1);

\draw[->] (R0S0R0) to (R0S0S0);
\draw[->] (R0S0R1) to (R0S0S0);

\draw[->] (R0R1R0) to (R0S0R0);
\draw[->] (R0R1R0) to (R0S0R1);

\draw[->] (R0R1R0) to (R0S1R0);
\draw[->] (R0R1R0) to (R0S1R1);

\draw[->] (R0S1R0) to (R0S1S0);
\draw[->] (R0S1R1) to (R0S1S0);

\draw[->] (R0R2R0) to (R0S1R0);
\draw[->] (R0R2R0) to (R0S1R1);

% Middle

\node (S0R0R0) at (6, 1) {$x$};

\node (S0S0R0) at (4, 2) {$x$};
\node (S0S0S0) at (5, 2) {$\alpha$};
\node (S0S0R1) at (6, 2) {$x$};
\node (S0S0S1) at (7, 2) {$\beta$};
\node (S0S0R2) at (8, 2) {$x$};

\node (S0R1R0) at (6, 3) {$x$};

\draw[->] (S0R0R0) to (S0S0R0);
\draw[->] (S0R0R0) to (S0S0R1);
\draw[->] (S0R0R0) to (S0S0R2);

\draw[->] (S0S0R0) to (S0S0S0);
\draw[->] (S0S0R1) to (S0S0S0);
\draw[->] (S0S0R1) to (S0S0S1);
\draw[->] (S0S0R2) to (S0S0S1);

\draw[->] (S0R1R0) to (S0S0R0);
\draw[->] (S0R1R0) to (S0S0R1);
\draw[->] (S0R1R0) to (S0S0R2);

% Target

\node (R1R0R0) at (11, 0) {$x$};

\node (R1S0R0) at (10, 1) {$x$};
\node (R1S0S0) at (11, 1) {$\beta$};
\node (R1S0R1) at (12, 1) {$x$};

\node (R1R1R0) at (11, 2) {$x$};

\node (R1S1R0) at (10, 3) {$x$};
\node (R1S1S0) at (11, 3) {$\alpha$};
\node (R1S1R1) at (12, 3) {$x$};

\node (R1R2R0) at (11, 4) {$x$};

\draw[->] (R1R0R0) to (R1S0R0);
\draw[->] (R1R0R0) to (R1S0R1);

\draw[->] (R1S0R0) to (R1S0S0);
\draw[->] (R1S0R1) to (R1S0S0);

\draw[->] (R1R1R0) to (R1S0R0);
\draw[->] (R1R1R0) to (R1S0R1);

\draw[->] (R1R1R0) to (R1S1R0);
\draw[->] (R1R1R0) to (R1S1R1);

\draw[->] (R1S1R0) to (R1S1S0);
\draw[->] (R1S1R1) to (R1S1S0);

\draw[->] (R1R2R0) to (R1S1R0);
\draw[->] (R1R2R0) to (R1S1R1);

% Slices

\node at (3, 2) {$\longrightarrow$};
\node at (9, 2) {$\longleftarrow$};

\end{tikzpicture}
\]
This is defined over the signature $\Sigma = \set{x, \alpha, \beta}$ with $\dim x = 0$, $\dim \alpha = 2$, and $\dim \beta = 2$.
\end{example}

\subsection{Explosion}

Given a diagram in the zigzag category $\Zig(\cat C)$, we can ``explode'' it to obtain a diagram in the base category $\cat C$.
This can be formally defined as the left adjoint of the zigzag construction; this result is inspired by Heidemann \cite{Hei23} and was first presented in our previous paper \cite{TV24}.

\begin{definition}
Given a category $\cat C$ with a functor $F : \cat C \to \Ord$, we define its \emph{explosion} $\Exp_F(\cat C)$ to the following category. 
\begin{itemize}
\item for every $x \in \cat C$ in the fibre of $n$, we have objects
\begin{align*}
&\sing^x_i &i &\in \ord n \\
&\reg^x_j &j &\in \ord{n + 1}
\end{align*}
\item for every $f : x \to y$ in the fibre of $\sigma : \ord n \to \ord m$ with $\rho = \Reg \sigma$, we have morphisms
\begin{align*}
\sing^x_i &\xrightarrow{\sing^f_i} \sing^y_{\sigma(i)} &i &\in \ord n \\
\reg^x_{\rho(j)} &\xrightarrow{\reg^f_j} \reg^y_j &j &\in \ord{m + 1} \\
\reg^x_j &\xrightarrow{\delta^f_{i, j}} \sing^y_i &i &\in \ord m, \, j \in \rho([i, i + 1])
\end{align*}
\end{itemize}
\end{definition}

\begin{proposition}
We have an adjunction
\[
\Exp : \Cat / \Ord \rightleftarrows \Cat : \Zig
\]
\end{proposition}

\begin{proof}
For every $\cat J \to \Ord$ and $\cat C$, we define a bijection
\begin{prooftree}
\AxiomC{$F : \cat J \to \Zig(\cat C)$}
\UnaryInfC{$\smallint F : \Exp(\cat J) \to \cat C$}
\end{prooftree}
Given a diagram $F$, we define the diagram $\smallint F$ as follows:
\begin{align*}
\smallint F(\sing^x_i) &\coloneqq F(x)(\sing i) \\
\smallint F(\reg^x_j) &\coloneqq F(x)(\reg j) \\
\smallint F(\sing^f_i) &\coloneqq F(f)(\sing i) \\
\smallint F(\reg^f_j) &\coloneqq F(f)(\reg j) \\
\smallint F(\delta^f_{i, j}) &\coloneqq F(f)(\reg j, \sing i)
\end{align*}
It is easy to check this is a natural isomorphism, so $\Exp \dashv \Zig$.
\end{proof}

By iterating the construction, we can also explode diagrams in the $n$-fold zigzag category:
\begin{prooftree}
\AxiomC{$F : \cat J \to \Zig^n(\cat C)$}
\UnaryInfC{$\smallint^k F : \Exp^k(\cat J) \to \Zig^{n - k}(\cat C)$}
\end{prooftree}
We call this the \emph{$k$-fold explosion}.
In particular, note that:
\begin{itemize}
\item If $X$ is an $n$-zigzag, seen as a functor $\ord 1 \to \Zig^n(\cat C)$, its $n$-fold explosion $\Exp^n(\ord 1) \to \cat C$ is $X$ seen as a diagram in $\cat C$.
\item If $f$ is an $n$-zigzag map, seen as a functor $\ord 2 \to \Zig^n(\cat C)$, its $n$-fold explosion $\Exp^n(\ord 2) \to \cat C$ is $f$ seen as a diagram in $\cat C$.
\end{itemize}

Given an $n$-zigzag $X$, we write $\Pt^k(X)$ for the $k$-fold explosion $\Exp^k(\ord 1)$.
The objects of $\Pt^k(X)$ are called the \emph{$k$-points} of $X$.
For convenience, we write these as $k$-tuples consisting of objects $\sing i$ and $\reg j$, i.e.\ we write
\begin{align*}
\sing^x_i &\leadsto (\sing i, x) \\
\reg^x_j &\leadsto (\reg j, x)
\end{align*}
In particular, the objects of $\Pt^k(X)$ can be inductively defined as follows:
\begin{align*}
\Pt^0(X) &= \set{ () } \\
\Pt^{k + 1}(X) &= \set{ (\sing i, x) \mid i \in X_\mathsf{s}, \, x \in \Pt^k(X(\sing i)) } \cup \set{ (\reg j, x) \mid j \in X_\mathsf{r}, \, x \in \Pt^k(X(\reg j)) }
\end{align*}
For example, the sequence of explosions of the monoid 2-diagram in \Cref{ex:monoid} is:
\[
()
\quad\leadsto\quad
\begin{aligned}\scalebox{1.0}{\begin{tikzpicture}

\node (R0) at (0, 0) {$(\reg 0)$};
\node (S0) at (0, 1) {$(\sing 0)$};
\node (R1) at (0, 2) {$(\reg 1)$};

\draw[->] (R0) to (S0);
\draw[->] (R1) to (S0);

\end{tikzpicture}}\end{aligned}
\quad\leadsto\quad
\begin{aligned}\scalebox{1.0}{\begin{tikzpicture}[xscale=1.75, yscale=1.5]

\node (R0R0) at (0, 0) {$(\reg 0, \reg 0)$};
\node (R0S0) at (1, 0) {$(\reg 0, \sing 0)$};
\node (R0R1) at (2, 0) {$(\reg 0, \reg 1)$};
\node (R0S1) at (3, 0) {$(\reg 0, \sing 1)$};
\node (R0R2) at (4, 0) {$(\reg 0, \reg 2)$};

\node (S0R0) at (0, 1) {$(\sing 0, \reg 0)$};
\node (S0S0) at (2, 1) {$(\sing 0, \sing 0)$};
\node (S0R1) at (4, 1) {$(\sing 0, \reg 1)$};

\node (R1R0) at (0, 2) {$(\reg 1, \reg 0)$};
\node (R1S0) at (2, 2) {$(\reg 1, \sing 0)$};
\node (R1R1) at (4, 2) {$(\reg 1, \reg 1)$};

\draw[->] (R0R0) to (R0S0);
\draw[->] (R0R1) to (R0S0);
\draw[->] (R0R1) to (R0S1);
\draw[->] (R0R2) to (R0S1);

\draw[->] (S0R0) to (S0S0);
\draw[->] (S0R1) to (S0S0);

\draw[->] (R1R0) to (R1S0);
\draw[->] (R1R1) to (R1S0);

\draw[->] (R0R0) to (S0R0);
\draw[->] (R0R2) to (S0R1);

\draw[->] (R0S0) to (S0S0);
\draw[->] (R0S1) to (S0S0);

\draw[->] (R1R0) to (S0R0);
\draw[->] (R1R1) to (S0R1);

\draw[->] (R1S0) to (S0S0);

\end{tikzpicture}}\end{aligned}
\]

Note that if $\cat J$ is a poset, then so is $\Exp_F(\cat J)$ for every $F : \cat J \to \Ord$.
Therefore, since $\ord 1$ is a poset, the category of $k$-points $\Pt^k(X)$ of an $n$-diagram $X$ is guaranteed to be a poset.

\section{Injectification} \label{sec:injectification}

In this section, we introduce a categorical construction, called \emph{injectification}, which can be thought of as a generalisation of factorisation.
We begin by recalling some definitions of factorisation systems~\cite{Rie08}.

\subsection{Factorization systems}

\begin{definition}
A \emph{functorial factorisation} on $\cat C$ is a functor
\[
\mathsf{fact}: [\ord 2, \cat C] \to [\ord 3, \cat C]
\]
that is a section of the composition functor $d_1: [\ord 3, \cat C] \to [\ord 2, \cat C]$.
\end{definition}

Explicitly, this specifies how to factorize morphisms in a way that depends functorially on commutative squares.
We write $\mathsf{fact} = (F, L, R)$ with the factorisations denoted as follows:
\[
\begin{tikzcd}
A \rar["f"] \dar & B \dar \\
C \rar["g"'] & D
\end{tikzcd}
\quad\mapsto\quad
\begin{tikzcd}
A \rar["L f"] \dar & F f \rar["R f"] \dar[dashed] & B \dar \\
C \rar["L g"'] & F g \rar["R g"'] & D
\end{tikzcd}
\]

Let $\mathcal{L}$ and $\mathcal{R}$ be two classes of morphisms of $\cat C$ closed under composition and containing all isomorphisms, corresponding to two wide subcategories $\cat C_\mathcal{L}$ and $\cat C_\mathcal{R}$ of $\cat C$.
We say that $(\mathcal{L}, \mathcal{R}, \mathsf{fact})$ forms a \emph{factorisation system} if $Lf \in \mathcal{L}$ and $Rf \in \mathcal{R}$ for every morphism $f$.

\begin{example}
\label{ex:mono-epi}
In $\Set$ and $\Pos$, there is a factorisation system given by:
\begin{itemize}
\item $\mathcal{L} = \text{monomorphisms}$ and $\mathcal{R} = \text{epimorphisms}$,
\item every map $f : A \to B$ is factored through the coproduct:
\[
\begin{tikzcd}
A \rar[>->] & A + B \rar[->>] & B
\end{tikzcd}
\]
\end{itemize}
\end{example}

\begin{example}
In any topos, there is a factorisation system given by:
\begin{itemize}
\item $\mathcal{L} = \text{epimorphisms}$ and $\mathcal{R} = \text{monomorphisms}$,
\item every morphism $f : A \to B$ is factored through its \emph{image}
\[
\begin{tikzcd}
A \rar[->>] & \operatorname{im} f \rar[>->] & B
\end{tikzcd}
\]
which is the smallest subobject of $B$ for which the factorisation exists.
\end{itemize}
\end{example}

\subsection{Main construction}

Fix a category $\cat C$ and a pair $(\mathcal{L}, \mathcal{R})$ of classes of morphisms of $\cat C$, both closed under isomorphisms.

\begin{definition}
\label{def:injectification}
Given a functor $X : \cat J \to \cat C$, an \emph{injectification} of $X$ is a functor
\[
\widehat{X} : \cat J \to \cat C_\mathcal{L}
\]
together with a natural transformation whose every component is in $\mathcal{R}$:
\[
\begin{tikzcd}
\cat J \ar[rr, "X", ""'{name=a, below}] \drar["\widehat{X}"'] && \cat C \\
& \cat C_\mathcal{L} \urar[hookrightarrow] \arrow[to=a, Rightarrow, "\epsilon"']
\end{tikzcd}
\]
\end{definition}

Conceptually, injectification takes a diagram and produces another diagram of the same shape, where every morphism has been replaced by one in $\mathcal{L}$, together with a way to map back to the original diagram by morphisms in $\mathcal{R}$.
The motivating example is when the two classes are (mono, epi), then injectification has the effect of making every morphism injective in a way that preserves the data of the original diagram.

We can also think of injectification as a way to generalise factorisations from morphisms to diagrams.
If we let $\cat J = \ord 2$, then a functor $\ord 2 \to \cat C$ picks out a morphism $f : A \to B$ in $\cat C$, and an injectification can be given by factorizing the morphism $f$ (if such a factorisation exists):
\[
\begin{tikzcd}
A \rar["f"] & B \\
A \uar[equal] \rar["L f"'] & F f \uar["R f"']
\end{tikzcd}
\]

The main theorem of this section shows that given a factorisation system, we can use it to construct injectifications of poset-shaped diagrams.
Moreover, the construction is functorial.

\begin{definition}
Let $J$ be a poset and let $F : J \to \cat C$ be a diagram such that $F i \to F j$ is in $\mathcal{L}$ for all $i \leq j$.
We say that $F$ is \emph{$\mathcal{L}$-faithful} if the colimit inclusions $F i \to \colim_J F$ are also in $\mathcal{L}$ for all $i \in J$, and $J$ is \emph{$\mathcal{L}$-faithful} if every diagram $F : J \to \cat C$ is $\mathcal{L}$-faithful.
\end{definition}

\begin{theorem}
\label{thm:injectification}
Let $J$ be a finite poset, and assume that:
\begin{itemize}
\item $\cat C$ is finitely cocomplete;
\item $(\mathcal{L}, \mathcal{R}, \mathsf{fact})$ is a factorisation system on $\cat C$ with $\mathsf{fact} = (F, L, R)$; and
\item every downward closed subset $I \subseteq J$ is $\mathcal{L}$-faithful.
\end{itemize}
Then any diagram $X : J \to \cat C$ has an injectification $(\widehat{X}, \epsilon)$.
Moreover, this construction is functorial, in that any natural transformation $\alpha : X \Rightarrow Y$ can be lifted to a commutative square in $[J, \cat C]$:
\[
\begin{tikzcd}
X \rar["\alpha"] & Y \\
\widehat{X} \uar["\epsilon_X"] \rar[dashed, "\widehat{\alpha}"'] & \widehat{Y} \uar["\epsilon_Y"']
\end{tikzcd}
\]
\end{theorem}

\begin{proof}
We construct the injectification by induction.\footnotemark{}
First, for any minimal $i \in J$, define 
\begin{align*}
\widehat{X}_i &\coloneqq X_i \\
\epsilon_i &\coloneqq \id_{X_i}
\end{align*}
For any non-minimal $i \in J$, restrict to the downward-closed subset
\[
J \downarrow i = \set{j \in J : j < i}.
\]
By induction, the injectification is already defined on $J \downarrow i$.
Hence, we have something like the following picture, depicting morphisms of $\mathcal{L}$ using $\rightarrowtail$, and morphisms of $\mathcal{R}$ using $\twoheadrightarrow$:
\[
\begin{tikzpicture}
\node (1) at (0, 0) {$X j_1$};
\node (2) at (2, 1) {$X j_2$};
\node (3) at (1, 2) {$X j_3$};

\node (0) at (5, 1) {$X i$};

\node (X1) at (0, -4) {$\widehat{X} j_1$};
\node (X2) at (2, -3) {$\widehat{X} j_2$};
\node (X3) at (1, -2) {$\widehat{X} j_3$};

\draw[->] (1) to (0);
\draw[->] (2) to (0);
\draw[->] (3) to (0);

\draw[->] (1) to (2);
\draw[->] (2) to (3);
\draw[->] (1) to (3);

\draw[->>] (X1) to (1);
\draw[->>] (X2) to (2);
\draw[->>] (X3) to (3);

\draw[>->] (X1) to (X2);
\draw[>->] (X2) to (X3);
\draw[>->] (X1) to (X3);

\node[draw,dotted,rounded corners,fit=(1) (2) (3)] {};
\node[draw,dotted,rounded corners,fit=(X1) (X2) (X3)] {};
\end{tikzpicture}
\]
Now take a colimit of $\widehat{X} : J \downarrow i \to \mathcal{L}$ in $\cat C$.
For every $j < i$, we have a morphism
\[
\widehat{X} j \xrightarrow{\epsilon_j} X j \to X i
\]
so by the universal property, there is a unique morphism $p : \colim_{J \downarrow i} \widehat{X} \to X i$, as follows:
\[
\begin{tikzpicture}
\node (1) at (0, 0) {$X j_1$};
\node (2) at (2, 1) {$X j_2$};
\node (3) at (1, 2) {$X j_3$};

\node (0) at (5, 1) {$X i$};

\node (X1) at (0, -4) {$\widehat{X} j_1$};
\node (X2) at (2, -3) {$\widehat{X} j_2$};
\node (X3) at (1, -2) {$\widehat{X} j_3$};

\node (colim) at (5, -3) {$\colim_{J \downarrow i} \widehat{X}$};

\draw[->] (1) to (0);
\draw[->] (2) to (0);
\draw[->] (3) to (0);

\draw[->] (1) to (2);
\draw[->] (2) to (3);
\draw[->] (1) to (3);

\draw[->>] (X1) to (1);
\draw[->>] (X2) to (2);
\draw[->>] (X3) to (3);

\draw[>->] (X1) to (X2);
\draw[>->] (X2) to (X3);
\draw[>->] (X1) to (X3);

\draw[>->] (X1) to (colim);
\draw[>->] (X2) to (colim);
\draw[>->] (X3) to (colim);

\draw[->, dashed] (colim) to (0);

\node[draw,dotted,rounded corners,fit=(1) (2) (3)] {};
\node[draw,dotted,rounded corners,fit=(X1) (X2) (X3)] {};
\end{tikzpicture}
\]
Now factorize the morphism $p$ as
\[
\colim_{J \downarrow i} \widehat{X} \xrightarrow{L p} F p \xrightarrow{R p} X i.
\]
Define $\widehat{X} i \coloneqq F p$.
Also for every $j < i$, define its image under $\widehat{X}$ to be the composite
\[
\widehat{X} j \to \colim_{J \downarrow i} \widehat{X} \xrightarrow{Lp} \widehat{X} i
\]
which is in $\mathcal{L}$ because the colimit inclusion is in $\mathcal{L}$.
Finally, define $\epsilon_i \coloneqq Rp$, and note that naturality follows trivially from the factorisation and the universal property of the colimit.
\[
\begin{tikzpicture}
\node (1) at (0, 0) {$X j_1$};
\node (2) at (2, 1) {$X j_2$};
\node (3) at (1, 2) {$X j_3$};

\node (0) at (8, 1) {$X i$};

\node (X1) at (0, -4) {$\widehat{X} j_1$};
\node (X2) at (2, -3) {$\widehat{X} j_2$};
\node (X3) at (1, -2) {$\widehat{X} j_3$};

\node (colim) at (5, -3) {$\colim_{J \downarrow i} \widehat{X}$};

\node (E) at (8, -3) {$\widehat{X} i$};

\draw[->] (1) to (0);
\draw[->] (2) to (0);
\draw[->] (3) to (0);

\draw[->] (1) to (2);
\draw[->] (2) to (3);
\draw[->] (1) to (3);

\draw[->>] (X1) to (1);
\draw[->>] (X2) to (2);
\draw[->>] (X3) to (3);

\draw[>->] (X1) to (X2);
\draw[>->] (X2) to (X3);
\draw[>->] (X1) to (X3);

\draw[>->] (X1) to (colim);
\draw[>->] (X2) to (colim);
\draw[>->] (X3) to (colim);

\draw[->, dashed] (colim) to (0);
\draw[>->] (colim) to (E);
\draw[->>] (E) to (0);

\node[draw,dotted,rounded corners,fit=(1) (2) (3)] {};
\node[draw,dotted,rounded corners,fit=(X1) (X2) (X3)] {};
\end{tikzpicture}
\]

To prove functoriality, we proceed by induction again.
For any minimal $i \in J$, let $\widehat{\alpha}_i \coloneqq \alpha_i$.
Now for any non-minimal $i \in J$, we have a natural transformation $\widehat{f} : \widehat{X} \Rightarrow \widehat{Y}$ on $J \downarrow i$ by induction.
By the universal property, this induces a morphism between the colimits such that the following diagram commutes:
\[
\begin{tikzcd}
\colim_{J \downarrow i} \widehat{X} \rar["p_X"] \dar & X i \dar["\alpha_i"] \\
\colim_{J \downarrow i} \widehat{Y} \rar["p_Y"'] & Y i
\end{tikzcd}
\]
Hence, since the factorisation is functorial, we get a morphism $\widehat{\alpha}_i$ with the required commutative square:
\[
\begin{tikzcd}
\widehat{X} i \rar["\epsilon_{X, i}"] \dar[dashed, "\widehat{\alpha}_i"'] & X i \dar["\alpha_i"] \\
\widehat{Y} i \rar["\epsilon_{Y, i}"'] & Y i
\end{tikzcd}
\]
\end{proof}

\footnotetext{This is an induction on the poset $J$. If $J$ is infinite, it is a transfinite induction.}

Note that faithfulness is a generalisation of the similarly named notion from Dokuchaev and Novikov's paper \cite{DN16};
their notion arises as the special case when $\cat C = \Set$ and $\mathcal{L}$ is the class of injective maps.
Not every poset is faithful.
For example, consider the following diagram in $\Set$ with $a \mapsto c$ and $b \mapsto d$:
\[
\begin{tikzcd}
\set{c, d} & \set{e} \\
\set{a} \uar \urar & \set{b} \uar \ular
\end{tikzcd}
\]
Every map in the diagram is injective, yet the colimit is the singleton set $\set{*}$ and the inclusion $\set{c, d} \to \set{*}$ is \emph{not} injective.
Therefore, the poset is not faithful with respect to injective maps.

Dokuchaev and Novikov proved that, for $\Set$, a poset is faithful (with respect to injections) iff it does not contain any \emph{subcrowns} satisfying a certain condition, where a crown is a poset of the form:
\[
\begin{tikzcd}[row sep=huge, column sep=small]
\bullet & \bullet & \cdots & \bullet \\
\bullet \uar \urar & \bullet \uar \urar & \cdots \urar & \bullet \uar \ar[ulll]
\end{tikzcd}
\]
Now the posets we are interested in are explosions $\Pt^k(X)$ of $n$-diagrams $X$, and it is easy to see that such posets never contain subcrowns, so they are faithful with respect to monomorphisms.

\section{Layout algorithm} \label{sec:algorithm}

We begin by defining the layout problem in more precise terms, using the explosion construction.

\begin{definition}
A \emph{full layout} for an $n$-diagram $X$ is a graph embedding of $\Pt^n(X)$ into~$\mathbb{R}^n$.
\end{definition}

\noindent
We are primarily interested in layouts that produce aesthetically pleasing string diagrams, since this forms an important part of \hio, where users must view and interact with terms graphically.
The goal of this section is to describe an algorithm that computes such layouts for arbitrary $n$-diagrams.

\subsection{Layout for zigzags}

First, let us consider the easier problem of layout for zigzags or 1-diagrams.

\begin{definition}
\label{def:layout}
A \emph{layout} for a zigzag $X$ consists of:
\begin{itemize}
\item a real number $w \in \mathbb{R}$ called the \emph{width}, together with
\item an injective monotone map $\lambda : \ord n \to [0, w]$ where $n$ is the length of $X$.
\end{itemize}
\end{definition}

This data can be extended to a map assigning positions to all objects of the zigzag, by interpolating the positions of the regular objects (except for the first and last). 
We define $\lambda' : \Pt(X) \to \mathbb{R}$ as follows:
\begin{equation}
\label{eq:layout}
\lambda'(\sing i) \coloneqq \lambda(i) + 1 \qquad
\lambda'(\reg j) \coloneqq \begin{cases*}
0 & \text{if $j = 0$} \\
\frac{1}{2}(\lambda(j - 1) + \lambda(j)) + 1 & \text{if $0 < j < n$} \\
w + 2 & \text{if $j = n$} \\
\end{cases*}
\end{equation}
If the zigzag has no singular objects, we set $\lambda'(\reg 0) \coloneqq \frac{1}{2} w + 1$ which is the midpoint of $[0, w + 2]$.

\begin{figure}
\centering
\begin{tikzpicture}
\node[lightgray] (R0) at (0, 0) {$\bullet$};
\node (S0) at (1, 0) {$0$};
\node[lightgray] (R1) at (2, 0) {$\bullet$};
\node (S1) at (3, 0) {$1$};
\node[lightgray] (R2) at (4, 0) {$\bullet$};

\draw[->, lightgray] (R0) to (S0);
\draw[->, lightgray] (R1) to (S0);
\draw[->, lightgray] (R1) to (S1);
\draw[->, lightgray] (R2) to (S1);

\node at (5, 0) {$\rightsquigarrow$};

\node (R0') at (6, 0) {$0$};
\node (S0') at (7, 0) {$1$};
\node (R1') at (8, 0) {$\frac{3}{2}$};
\node (S1') at (9, 0) {$2$};
\node (R2') at (10, 0) {$3$};

\draw[->] (R0') to (S0');
\draw[->] (R1') to (S0');
\draw[->] (R1') to (S1');
\draw[->] (R2') to (S1');
\end{tikzpicture}
\caption{A layout for a zigzag.}
\label{fig:layout-1d}
\end{figure}

More generally, we will want to compute a layout for a family of zigzags subject to certain aesthetic conditions relating the layouts of each zigzag according to some zigzag maps between them.

\begin{definition}
\label{def:fair-layout}
Given a functor $X : \cat J \to \Zig(\cat C)$, a \emph{layout} consists of:
\begin{itemize}
\item a real number $w \in \mathbb{R}$ called the \emph{width}, together with
\item a family of injective monotone maps $\lambda_i : (X i)_\mathsf{s} \to [0, w]$ for every $i \in \cat J$. 
\end{itemize}
Given a (non-identity) morphism $f : i \to j$ in $\cat J$, we say that the layout is \emph{$f$-fair} if
\[
\forall q \in (X j)_\mathsf{s} \ldotp \lambda_j(q) \text{ is the average of } \set{ \lambda_i(p) \mid p \in (X i)_\mathsf{s} \text { such that } (X f)_\mathsf{s}(p) = q }.
\]
If a layout is fair for every (non-identity) morphism, we simply say that the layout is \emph{fair}.
\end{definition}

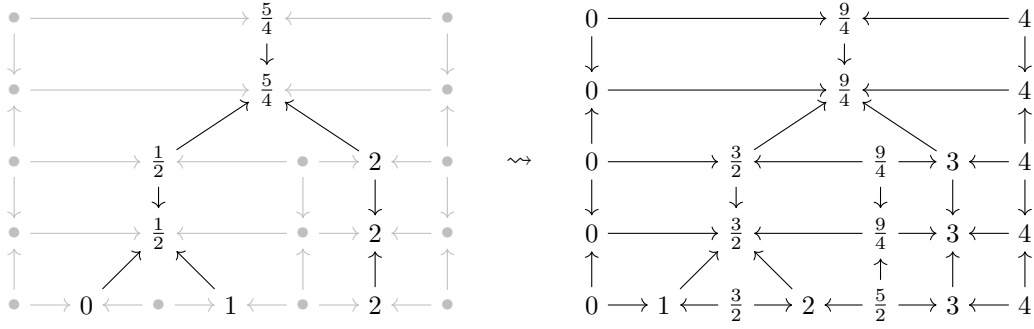
\begin{figure}
\centering
\begin{tikzpicture}[scale=0.95]
\node[lightgray] (R0R0) at (0, 0) {$\bullet$};
\node (R0S0) at (1, 0) {$0$};
\node[lightgray] (R0R1) at (2, 0) {$\bullet$};
\node (R0S1) at (3, 0) {$1$};
\node[lightgray] (R0R2) at (4, 0) {$\bullet$};
\node (R0S2) at (5, 0) {$2$};
\node[lightgray] (R0R3) at (6, 0) {$\bullet$};

\node[lightgray] (S0R0) at (0, 1) {$\bullet$};
\node (S0S0) at (2, 1) {$\frac{1}{2}$};
\node[lightgray] (S0R1) at (4, 1) {$\bullet$};
\node (S0S1) at (5, 1) {$2$};
\node[lightgray] (S0R2) at (6, 1) {$\bullet$};

\node[lightgray] (R1R0) at (0, 2) {$\bullet$};
\node (R1S0) at (2, 2) {$\frac{1}{2}$};
\node[lightgray] (R1R1) at (4, 2) {$\bullet$};
\node (R1S1) at (5, 2) {$2$};
\node[lightgray] (R1R2) at (6, 2) {$\bullet$};

\node[lightgray] (S1R0) at (0, 3) {$\bullet$};
\node (S1S0) at (3.5, 3) {$\frac{5}{4}$};
\node[lightgray] (S1R1) at (6, 3) {$\bullet$};

\node[lightgray] (R2R0) at (0, 4) {$\bullet$};
\node (R2S0) at (3.5, 4) {$\frac{5}{4}$};
\node[lightgray] (R2R1) at (6, 4) {$\bullet$};

\draw[->, lightgray] (R0R0) to (R0S0);
\draw[->, lightgray] (R0R1) to (R0S0);
\draw[->, lightgray] (R0R1) to (R0S1);
\draw[->, lightgray] (R0R2) to (R0S1);
\draw[->, lightgray] (R0R2) to (R0S2);
\draw[->, lightgray] (R0R3) to (R0S2);

\draw[->, lightgray] (S0R0) to (S0S0);
\draw[->, lightgray] (S0R1) to (S0S0);
\draw[->, lightgray] (S0R1) to (S0S1);
\draw[->, lightgray] (S0R2) to (S0S1);

\draw[->, lightgray] (R1R0) to (R1S0);
\draw[->, lightgray] (R1R1) to (R1S0);
\draw[->, lightgray] (R1R1) to (R1S1);
\draw[->, lightgray] (R1R2) to (R1S1);

\draw[->, lightgray] (S1R0) to (S1S0);
\draw[->, lightgray] (S1R1) to (S1S0);

\draw[->, lightgray] (R2R0) to (R2S0);
\draw[->, lightgray] (R2R1) to (R2S0);

\draw[->, lightgray] (R0R0) to (S0R0);
\draw[->, lightgray] (R0R2) to (S0R1);
\draw[->, lightgray] (R0R3) to (S0R2);

\draw[->] (R0S0) to (S0S0);
\draw[->] (R0S1) to (S0S0);
\draw[->] (R0S2) to (S0S1);

\draw[->, lightgray] (R1R0) to (S0R0);
\draw[->, lightgray] (R1R1) to (S0R1);
\draw[->, lightgray] (R1R2) to (S0R2);

\draw[->] (R1S0) to (S0S0);
\draw[->] (R1S1) to (S0S1);

\draw[->, lightgray] (R1R0) to (S1R0);
\draw[->, lightgray] (R1R2) to (S1R1);

\draw[->] (R1S0) to (S1S0);
\draw[->] (R1S1) to (S1S0);

\draw[->, lightgray] (R2R0) to (S1R0);
\draw[->, lightgray] (R2R1) to (S1R1);

\draw[->] (R2S0) to (S1S0);

\node at (7, 2) {$\rightsquigarrow$};

\node (R0R0') at (8, 0) {$0$};
\node (R0S0') at (9, 0) {$1$};
\node (R0R1') at (10, 0) {$\frac{3}{2}$};
\node (R0S1') at (11, 0) {$2$};
\node (R0R2') at (12, 0) {$\frac{5}{2}$};
\node (R0S2') at (13, 0) {$3$};
\node (R0R3') at (14, 0) {$4$};

\node (S0R0') at (8, 1) {$0$};
\node (S0S0') at (10, 1) {$\frac{3}{2}$};
\node (S0R1') at (12, 1) {$\frac{9}{4}$};
\node (S0S1') at (13, 1) {$3$};
\node (S0R2') at (14, 1) {$4$};

\node (R1R0') at (8, 2) {$0$};
\node (R1S0') at (10, 2) {$\frac{3}{2}$};
\node (R1R1') at (12, 2) {$\frac{9}{4}$};
\node (R1S1') at (13, 2) {$3$};
\node (R1R2') at (14, 2) {$4$};

\node (S1R0') at (8, 3) {$0$};
\node (S1S0') at (11.5, 3) {$\frac{9}{4}$};
\node (S1R1') at (14, 3) {$4$};

\node (R2R0') at (8, 4) {$0$};
\node (R2S0') at (11.5, 4) {$\frac{9}{4}$};
\node (R2R1') at (14, 4) {$4$};

\draw[->] (R0R0') to (R0S0');
\draw[->] (R0R1') to (R0S0');
\draw[->] (R0R1') to (R0S1');
\draw[->] (R0R2') to (R0S1');
\draw[->] (R0R2') to (R0S2');
\draw[->] (R0R3') to (R0S2');

\draw[->] (S0R0') to (S0S0');
\draw[->] (S0R1') to (S0S0');
\draw[->] (S0R1') to (S0S1');
\draw[->] (S0R2') to (S0S1');

\draw[->] (R1R0') to (R1S0');
\draw[->] (R1R1') to (R1S0');
\draw[->] (R1R1') to (R1S1');
\draw[->] (R1R2') to (R1S1');

\draw[->] (S1R0') to (S1S0');
\draw[->] (S1R1') to (S1S0');

\draw[->] (R2R0') to (R2S0');
\draw[->] (R2R1') to (R2S0');

\draw[->] (R0R0') to (S0R0');
\draw[->] (R0R2') to (S0R1');
\draw[->] (R0R3') to (S0R2');

\draw[->] (R0S0') to (S0S0');
\draw[->] (R0S1') to (S0S0');
\draw[->] (R0S2') to (S0S1');

\draw[->] (R1R0') to (S0R0');
\draw[->] (R1R1') to (S0R1');
\draw[->] (R1R2') to (S0R2');

\draw[->] (R1S0') to (S0S0');
\draw[->] (R1S1') to (S0S1');

\draw[->] (R1R0') to (S1R0');
\draw[->] (R1R2') to (S1R1');

\draw[->] (R1S0') to (S1S0');
\draw[->] (R1S1') to (S1S0');

\draw[->] (R2R0') to (S1R0');
\draw[->] (R2R1') to (S1R1');

\draw[->] (R2S0') to (S1S0');
\end{tikzpicture}
\caption{A fair layout for an exploded 2-diagram.}
\label{fig:layout-2d}
\end{figure}

We are primarily interested in functors of the form $\Pt^n(X) \to \Zig(\cat C)$ which are the $n$-fold explosion of an $(n + 1)$-diagram $X$.
A layout for such a functor can be extended, using \eqref{eq:layout}, to a map $\Pt^{n + 1}(X) \to \mathbb{R}$.
While this does not give a full layout for $X$, it is like assigning $x$-coordinates to every point of $X$ and is a crucial step in computing a full layout.

The fairness condition ensures that points are centred with respect to their inputs and outputs, and in particular, that identities are rendered as straight wires or flat surfaces.
In general, it is not possible to get a fair layout for $\Pt^n(X) \to \Zig(\cat C)$ for $n \geq 2$.
For example, in the associator 3-diagram in \Cref{fig:associator}, it is not possible to assign $x$-coordinates in a way that makes the top surface flat in both the $y$- and $z$-axes.
However, what we are looking for is a layout that is as close to fair as possible, which becomes an optimisation problem.

\paragraph*{Zigzag layout algorithm}
Given a finite poset-shaped diagram
\[
X : J \to \Zig(\cat C)
\]
we obtain a (fair) layout according to the following scheme:
\begin{enumerate}[(1)]
\item Build the diagram $J \xrightarrow{X} \Zig(\cat C) \to \Ord \hookrightarrow \Pos$ which we denote by $P$.
\item Obtain an injectification using \Cref{thm:injectification} for the factorisation system in \Cref{ex:mono-epi}:
\[
\begin{tikzcd}
J \ar[rr, "P", ""'{name=a, below}] \drar["\widehat{P}"'] && \Pos \\
& \Pos_\mathrm{inj} \urar[hookrightarrow] \arrow[to=a, Rightarrow, "\epsilon"']
\end{tikzcd}
\]
Assume wlog that $P i \subseteq \widehat{P} i$ for all $i \in J$ and $\widehat{P} i \subseteq \widehat{P} j$ whenever $i < j$. 
\item Construct a linear programming problem using the following scheme:\footnotemark{}
\begin{align*}
\textbf{minimise} \quad
& w + \sum m_{i, j, p} \\
\textbf{subject to} \quad
& 0 \leq v_{i, p} \leq w & ( i \in J, \, p \in \widehat{P} i ) \quad (\star) \\
& v_{i, q} - v_{i, p} \geq 1 & ( p < q \text{ in } \widehat{P} i ) \quad (\ast) \\
& v_{i, q} - v_{i, p} \geq 1 & ( \epsilon_i(p) < \epsilon_j(q) \text{ in } P i ) \quad (\diamond) \\
& \abs{v_{i, p} - v_{j, p}} \leq m_{i, j, p} & ( i < j ) \quad (\dag) \\
& v_{j, q} = \operatorname{avg}(\set{v_{j, p} \mid p \in P i, \, P_{i < j}(p) = q}) & ( i < j, \, q \in P j ) \quad (\ddag)
\end{align*}
\item Extract the layout from the solution: $\lambda_i(p) \coloneqq v_{i, p}$ for $i \in J$ and $p \in P i$.
\end{enumerate}

\footnotetext{In the implementation, we solve this using HiGHS, an open-source linear programming solver~\cite{HH18}.}

The idea is as follows.
The set $\widehat{P} i$ is a colimit union of all $P j$ for $j \leq i$.
Therefore, we have a variable $v_{i, p}$ for every $p \in P j$ for $j \leq i$.
We think of this as the value of $\lambda_j(p)$ locally at $i$, with the true value being the one at $j$.
This gives for each $i$ a \emph{local} layout for the points in the neighbourhood of $i$. The constraints $(\dag)$ assert that the local layouts agree on common points, whereas the constraints $(\ddag)$ impose that each local layout is fair, and together, they imply that the global layout is fair.
However, this may not in general be satisfiable, so the constraints $(\dag)$ are weakened so that the difference between local layouts is minimised instead of made zero.
This gives a layout that, while not fair, is as fair as possible.

Injectification plays a key role in computing these constraints.
For this to work, we had to move from $\Ord$ to $\Pos$, the category of finite posets, so that we are working in a category with sufficient colimit structure.
Recently, a notion of zigzags natively indexed over $\Pos$ instead of $\Ord$ was given by Sarti and Vicary \cite{SV23}, and it would be interesting to see how the layout algorithm can be generalised to this case.

\subsection{Layout for diagrams}

We can obtain a full layout for any $n$-diagram $X$ recursively:
\begin{enumerate}[(1)]
\item If $X$ is a 0-diagram, its layout is trivial:
\begin{align*}
\lambda: \Pt^0(X) &\to \mathbb{R}^0 \\
() &\mapsto ()
\end{align*}
\item If $X$ is an $(n + 1)$-diagram, compute its layout as an $n$-diagram:
\[
\lambda_\mathrm{ind} : \Pt^n(X) \to \mathbb{R}^n
\]
\item Construct its $n$-fold explosion $\Pt^n(X) \to \Zig(\cat C)$.
\item Obtain its layout using the zigzag layout algorithm:
\[
\lambda_\mathrm{zig} : \Pt^{n + 1}(X) \to \mathbb{R}
\]
\item Finally, combine this with the inductive layout to get a full layout for $X$:
\begin{align*}
\lambda : \Pt^{n + 1}(X) &\to \mathbb{R}^{n + 1} \\
(\sing i, x) &\mapsto (\lambda_\mathrm{zig}(\sing i, x), \lambda_\mathrm{ind}(x)) \\
(\reg j, x) &\mapsto (\lambda_\mathrm{zig}(\reg j, x), \lambda_\mathrm{ind}(x))
\end{align*}
\end{enumerate}

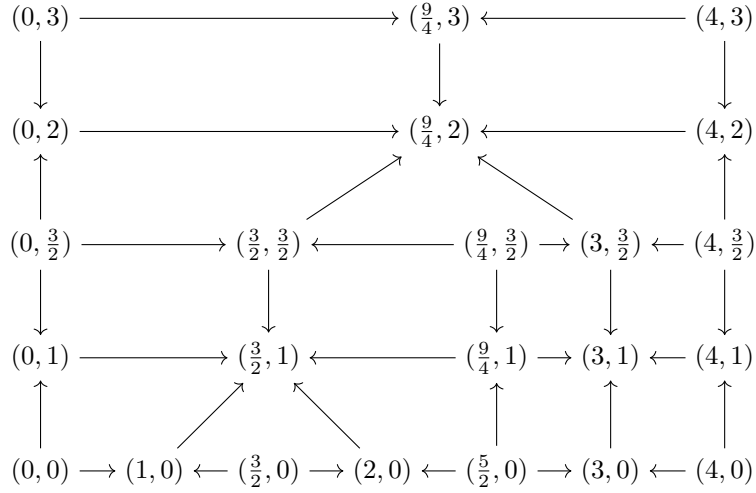
\begin{figure}
\centering
\begin{tikzpicture}[scale=1.5]
\node (R0R0) at (0, 0) {$(0, 0)$};
\node (R0S0) at (1, 0) {$(1, 0)$};
\node (R0R1) at (2, 0) {$(\frac{3}{2}, 0)$};
\node (R0S1) at (3, 0) {$(2, 0)$};
\node (R0R2) at (4, 0) {$(\frac{5}{2}, 0)$};
\node (R0S2) at (5, 0) {$(3, 0)$};
\node (R0R3) at (6, 0) {$(4, 0)$};

\node (S0R0) at (0, 1) {$(0, 1)$};
\node (S0S0) at (2, 1) {$(\frac{3}{2}, 1)$};
\node (S0R1) at (4, 1) {$(\frac{9}{4}, 1)$};
\node (S0S1) at (5, 1) {$(3, 1)$};
\node (S0R2) at (6, 1) {$(4, 1)$};

\node (R1R0) at (0, 2) {$(0, \frac{3}{2})$};
\node (R1S0) at (2, 2) {$(\frac{3}{2}, \frac{3}{2})$};
\node (R1R1) at (4, 2) {$(\frac{9}{4}, \frac{3}{2})$};
\node (R1S1) at (5, 2) {$(3, \frac{3}{2})$};
\node (R1R2) at (6, 2) {$(4, \frac{3}{2})$};

\node (S1R0) at (0, 3) {$(0, 2)$};
\node (S1S0) at (3.5, 3) {$(\frac{9}{4}, 2)$};
\node (S1R1) at (6, 3) {$(4, 2)$};

\node (R2R0) at (0, 4) {$(0, 3)$};
\node (R2S0) at (3.5, 4) {$(\frac{9}{4}, 3)$};
\node (R2R1) at (6, 4) {$(4, 3)$};

\draw[->] (R0R0) to (R0S0);
\draw[->] (R0R1) to (R0S0);
\draw[->] (R0R1) to (R0S1);
\draw[->] (R0R2) to (R0S1);
\draw[->] (R0R2) to (R0S2);
\draw[->] (R0R3) to (R0S2);

\draw[->] (S0R0) to (S0S0);
\draw[->] (S0R1) to (S0S0);
\draw[->] (S0R1) to (S0S1);
\draw[->] (S0R2) to (S0S1);

\draw[->] (R1R0) to (R1S0);
\draw[->] (R1R1) to (R1S0);
\draw[->] (R1R1) to (R1S1);
\draw[->] (R1R2) to (R1S1);

\draw[->] (S1R0) to (S1S0);
\draw[->] (S1R1) to (S1S0);

\draw[->] (R2R0) to (R2S0);
\draw[->] (R2R1) to (R2S0);

\draw[->] (R0R0) to (S0R0);
\draw[->] (R0R2) to (S0R1);
\draw[->] (R0R3) to (S0R2);

\draw[->] (R0S0) to (S0S0);
\draw[->] (R0S1) to (S0S0);
\draw[->] (R0S2) to (S0S1);

\draw[->] (R1R0) to (S0R0);
\draw[->] (R1R1) to (S0R1);
\draw[->] (R1R2) to (S0R2);

\draw[->] (R1S0) to (S0S0);
\draw[->] (R1S1) to (S0S1);

\draw[->] (R1R0) to (S1R0);
\draw[->] (R1R2) to (S1R1);

\draw[->] (R1S0) to (S1S0);
\draw[->] (R1S1) to (S1S0);

\draw[->] (R2R0) to (S1R0);
\draw[->] (R2R1) to (S1R1);

\draw[->] (R2S0) to (S1S0);
\end{tikzpicture}
\caption{A full layout for a 2-diagram (combining the layouts of \Cref{fig:layout-1d,fig:layout-2d}).}
\end{figure}

\bibliographystyle{plainurl}
\bibliography{Layout}

\newpage
\appendix

\section{Further examples} \label{sec:examples}

Here we present more examples of 3D and 4D diagrams rendered using our layout algorithm.

\subsection{Figure-eight isotopy (\href{https://beta.homotopy.io/p/2402.00003}{link to workspace})}

Consider the associative 4-category generated by:
\begin{itemize}
\item a 0-cell $x$, and
\item an invertible 2-cell $f : \id x \to \id x$.
\end{itemize}
We can construct a homotopy 4-cell between two `figure-eight' 3-cells:
\[
\begin{aligned}\scalebox{0.4}{\begin{tikzpicture}
\definecolor{generator-1-2-1-zer}{RGB}{235, 170, 164}
\definecolor{generator-1-2-0-neg}{RGB}{121, 36, 27}
\definecolor{generator-1-2-0-pos}{RGB}{192, 57, 43}
\definecolor{generator-0-0-1-pos}{RGB}{56, 150, 211}

\newcommand{\wire}[2]{
  \ifdefined\recolor\draw[color=\recolor, line width=10pt]\else\draw[color=#1, line width=5pt]\fi #2
}
\newcommand{\clipped}[3]{
\begin{scope}
  \newcommand{\recolor}{#1}
  \clip#3;
  #2
\end{scope}
}

\begin{scope}[transparency group]
% Background surfaces
\fill[generator-0-0-1-pos] (0,0) -- (6,0) -- (6,8) -- (0,8) -- (0,0);
\newcommand{\layer}[1]{
  \clipped{generator-0-0-1-pos}{#1}{(0,0) -- (6,0) -- (6,8) -- (0,8) -- (0,0)}
  #1
}

% Wire layers
\wire{generator-1-2-0-neg}{(2,3) .. controls (2,3.8) and (2.4,4) .. (3,4) .. controls (3.6,4) and (4,4.2) .. (4,5)};
\layer{
\wire{generator-1-2-0-pos}{(3,2) .. controls (3.6,2) and (4,2.2) .. (4,3) .. controls (4,3.8) and (3.6,4) .. (3,4) .. controls (2.4,4) and (2,4.2) .. (2,5) .. controls (2,5.8) and (2.4,6) .. (3,6)};
\wire{generator-1-2-0-neg}{(3,2) .. controls (2.4,2) and (2,2.2) .. (2,3)(4,5) .. controls (4,5.8) and (3.6,6) .. (3,6)};
}
\end{scope}
\fill[generator-1-2-1-zer] (3,2) circle (0.14);
\fill[generator-1-2-1-zer] (3,6) circle (0.14);
\end{tikzpicture}
}\end{aligned}
\quad\to\quad
\begin{aligned}\scalebox{0.4}{\begin{tikzpicture}
\definecolor{generator-1-2-1-zer}{RGB}{235, 170, 164}
\definecolor{generator-1-2-0-pos}{RGB}{192, 57, 43}
\definecolor{generator-1-2-0-neg}{RGB}{121, 36, 27}
\definecolor{generator-0-0-1-pos}{RGB}{56, 150, 211}

\newcommand{\wire}[2]{
  \ifdefined\recolor\draw[color=\recolor, line width=10pt]\else\draw[color=#1, line width=5pt]\fi #2
}
\newcommand{\clipped}[3]{
\begin{scope}
  \newcommand{\recolor}{#1}
  \clip#3;
  #2
\end{scope}
}

\begin{scope}[transparency group]
% Background surfaces
\fill[generator-0-0-1-pos] (0,0) -- (6,0) -- (6,8) -- (0,8) -- (0,0);
\newcommand{\layer}[1]{
  \clipped{generator-0-0-1-pos}{#1}{(0,0) -- (6,0) -- (6,8) -- (0,8) -- (0,0)}
  #1
}

% Wire layers
\wire{generator-1-2-0-pos}{(2,3) .. controls (2,3.8) and (2.4,4) .. (3,4) .. controls (3.6,4) and (4,4.2) .. (4,5)};
\layer{
\wire{generator-1-2-0-pos}{(3,2) .. controls (2.4,2) and (2,2.2) .. (2,3)(4,5) .. controls (4,5.8) and (3.6,6) .. (3,6)};
\wire{generator-1-2-0-neg}{(3,2) .. controls (3.6,2) and (4,2.2) .. (4,3) .. controls (4,3.8) and (3.6,4) .. (3,4) .. controls (2.4,4) and (2,4.2) .. (2,5) .. controls (2,5.8) and (2.4,6) .. (3,6)};
}
\end{scope}
\fill[generator-1-2-1-zer] (3,2) circle (0.14);
\fill[generator-1-2-1-zer] (3,6) circle (0.14);
\end{tikzpicture}
}\end{aligned}
\]
This has been rendered in 4D as a movie of 3D diagrams, and can be viewed on YouTube:
\begin{center}
\url{https://youtu.be/n0K2TT7I0oM}
\end{center}

\subsection{Figure-eight is self-invertible (\href{https://beta.homotopy.io/p/2402.00004}{link to workspace})}

Consider the associative 5-category generated by:
\begin{itemize}
\item a 0-cell $x$, and
\item an invertible 3-cell $f : \id(\id x) \to \id(\id x)$.
\end{itemize}
We can construct a homotopy 5-cell showing that the `figure-eight' 4-cell is self-invertible:
\[
\begin{aligned}\scalebox{0.4}{\begin{tikzpicture}
\definecolor{generator-1-3-0-pos}{RGB}{192, 57, 43}
\definecolor{generator-1-3-0-neg}{RGB}{121, 36, 27}
\definecolor{generator-0-0-2-pos}{RGB}{89, 167, 218}
\definecolor{generator-1-3-1-zer}{RGB}{235, 170, 164}

\newcommand{\wire}[2]{
  \ifdefined\recolor\draw[color=\recolor, line width=10pt]\else\draw[color=#1, line width=5pt]\fi #2
}
\newcommand{\clipped}[3]{
\begin{scope}
  \newcommand{\recolor}{#1}
  \clip#3;
  #2
\end{scope}
}

\begin{scope}[transparency group]
% Background surfaces
\fill[generator-0-0-2-pos] (0,0) -- (6,0) -- (6,14) -- (0,14) -- (0,0);
\newcommand{\layer}[1]{
  \clipped{generator-0-0-2-pos}{#1}{(0,0) -- (6,0) -- (6,14) -- (0,14) -- (0,0)}
  #1
}

% Wire layers
\wire{generator-1-3-0-neg}{(2,3) .. controls (2,3.8) and (2.4,4) .. (3,4) .. controls (3.6,4) and (4,4.2) .. (4,5)(2,9) .. controls (2,9.8) and (2.4,10) .. (3,10) .. controls (3.6,10) and (4,10.2) .. (4,11)};
\layer{
\wire{generator-1-3-0-pos}{(3,2) .. controls (3.6,2) and (4,2.2) .. (4,3) .. controls (4,3.8) and (3.6,4) .. (3,4) .. controls (2.4,4) and (2,4.2) .. (2,5) .. controls (2,5.8) and (2.4,6) .. (3,6)(3,8) .. controls (3.6,8) and (4,8.2) .. (4,9) .. controls (4,9.8) and (3.6,10) .. (3,10) .. controls (2.4,10) and (2,10.2) .. (2,11) .. controls (2,11.8) and (2.4,12) .. (3,12)};
\wire{generator-1-3-0-neg}{(3,2) .. controls (2.4,2) and (2,2.2) .. (2,3)(4,5) .. controls (4,5.8) and (3.6,6) .. (3,6)(3,8) .. controls (2.4,8) and (2,8.2) .. (2,9)(4,11) .. controls (4,11.8) and (3.6,12) .. (3,12)};
}
\end{scope}
\fill[generator-1-3-1-zer] (3,2) circle (0.14);
\fill[generator-1-3-1-zer] (3,6) circle (0.14);
\fill[generator-1-3-1-zer] (3,8) circle (0.14);
\fill[generator-1-3-1-zer] (3,12) circle (0.14);
\end{tikzpicture}
}\end{aligned}
\quad\to\quad
\begin{aligned}\scalebox{0.4}{\begin{tikzpicture}
\definecolor{generator-0-0-2-pos}{RGB}{89, 167, 218}
\begin{scope}
% Background surfaces
\fill[generator-0-0-2-pos] (0,0) -- (2,0) -- (2,2) -- (0,2) -- (0,0);
% Wire layers
\end{scope}
\end{tikzpicture}
}\end{aligned}
\]
This has been rendered in 4D as a movie of 3D diagrams, and can be viewed on YouTube:
\begin{center}
\url{https://youtu.be/TWvfrXa-eQ0}
\end{center}

\subsection{Swallowtail coherence (\href{https://beta.homotopy.io/p/2402.00005}{link to workspace})}

\begin{figure}
\centering
\includegraphics[width=0.4\textwidth]{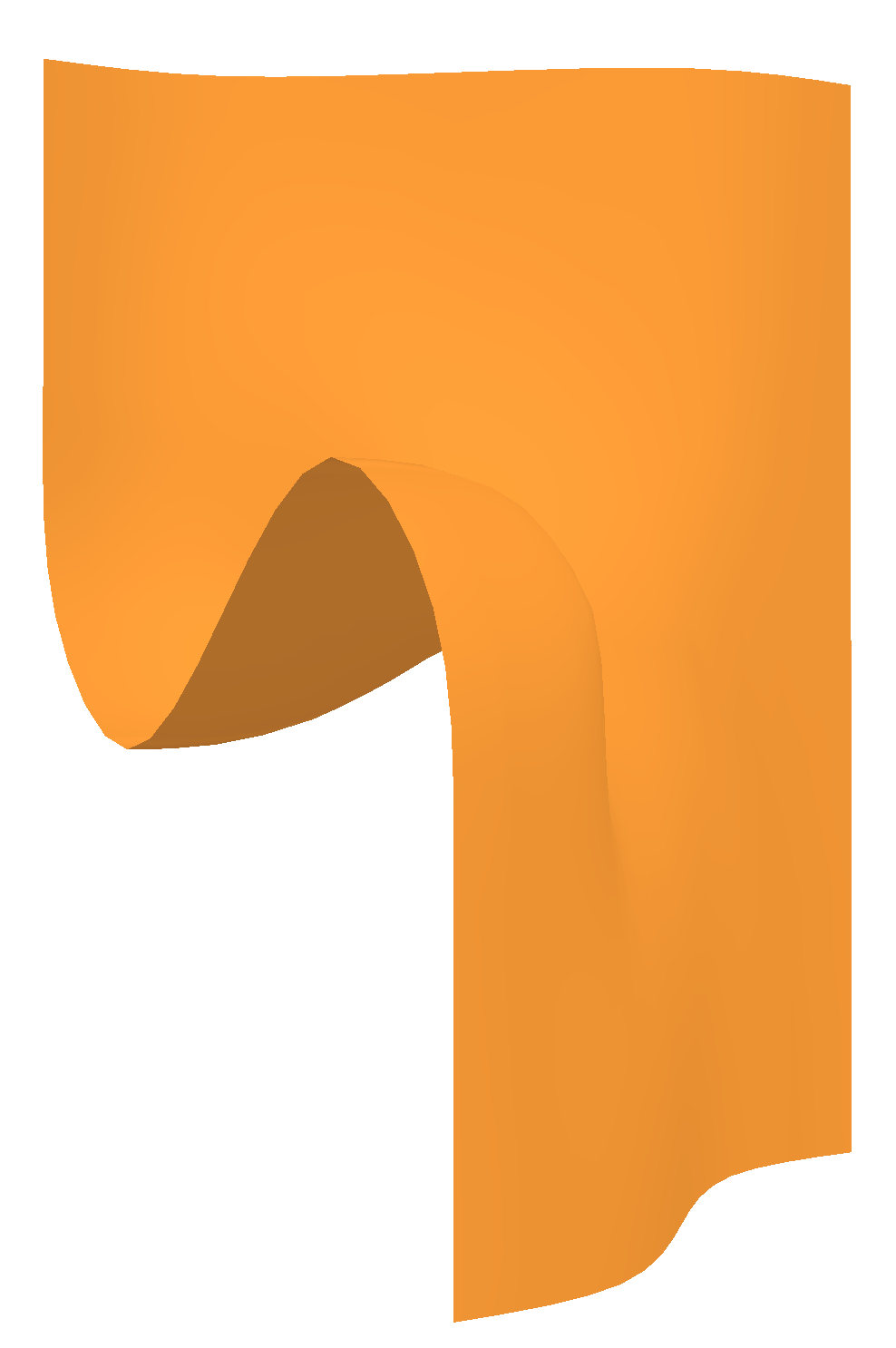}
\caption{The cusp rendered in 3D.}
\label{fig:cusp-3d}
\end{figure}

\begin{figure}
\centering
\includegraphics[width=0.5\textwidth]{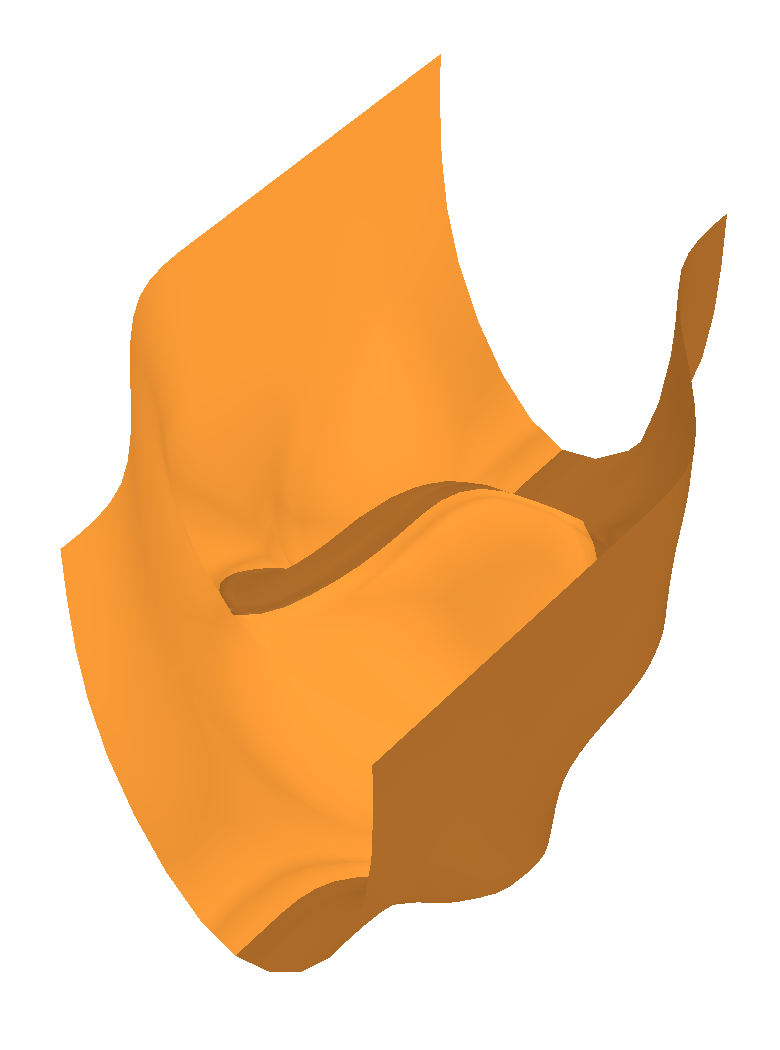}
\caption{The source of the swallowtail rendered in 3D.}
\label{fig:swallowtail-source-3d}
\end{figure}

Consider the associative 4-category generated by:
\begin{itemize}
\item a pair of 0-cells $x$ and $y$, and
\item an invertible 1-cell $f : x \to y$.
\end{itemize}
In particular, the invertibility implies that we have the following data:\footnotemark{}
\begin{itemize}
\item a 1-cell $f^{-1} : y \to y$, and
\item invertible 2-cells known as \emph{unit} and \emph{counit}:
\[
\begin{aligned}\scalebox{0.4}{\begin{tikzpicture}
\definecolor{generator-1-0-0-pos}{RGB}{192, 57, 43}
\definecolor{generator-2-1-0-neg}{RGB}{166, 105, 8}
\definecolor{generator-2-1-0-pos}{RGB}{243, 156, 18}
\definecolor{generator-0-0-0-pos}{RGB}{41, 128, 185}
\definecolor{generator-2-1-1-zer}{RGB}{251, 221, 173}
\begin{scope}
% Background surfaces
\fill[generator-0-0-0-pos] (0,0) -- (2,0) -- (2,1) .. controls (2,1.8) and (2.4,2) .. (3,2) .. controls (3.6,2) and (4,1.8) .. (4,1) -- (4,0) -- (6,0) -- (6,4) -- (0,4) -- (0,0);
\fill[generator-1-0-0-pos] (2,0) -- (4,0) -- (4,1) .. controls (4,1.8) and (3.6,2) .. (3,2) .. controls (2.4,2) and (2,1.8) .. (2,1) -- (2,0);
% Wire layers
\draw[color=generator-2-1-0-neg, line width=5pt](4,0) -- (4,1) .. controls (4,1.8) and (3.6,2) .. (3,2);
\draw[color=generator-2-1-0-pos, line width=5pt](2,0) -- (2,1) .. controls (2,1.8) and (2.4,2) .. (3,2);
\end{scope}
\fill[generator-2-1-1-zer] (3,2) circle (0.14);
\end{tikzpicture}}\end{aligned}
\]
\[
\begin{aligned}\scalebox{0.4}{\begin{tikzpicture}
\definecolor{generator-1-0-0-pos}{RGB}{192, 57, 43}
\definecolor{generator-2-1-0-neg}{RGB}{166, 105, 8}
\definecolor{generator-2-1-0-pos}{RGB}{243, 156, 18}
\definecolor{generator-0-0-0-pos}{RGB}{41, 128, 185}
\definecolor{generator-2-1-1-zer}{RGB}{251, 221, 173}
\begin{scope}
% Background surfaces
\fill[generator-0-0-0-pos] (2,3) .. controls (2,2.2) and (2.4,2) .. (3,2) .. controls (3.6,2) and (4,2.2) .. (4,3) -- (4,4) -- (2,4) -- (2,3);
\fill[generator-1-0-0-pos] (0,0) -- (6,0) -- (6,4) -- (4,4) -- (4,3) .. controls (4,2.2) and (3.6,2) .. (3,2) .. controls (2.4,2) and (2,2.2) .. (2,3) -- (2,4) -- (0,4) -- (0,0);
% Wire layers
\draw[color=generator-2-1-0-neg, line width=5pt](3,2) .. controls (2.4,2) and (2,2.2) .. (2,3) -- (2,4);
\draw[color=generator-2-1-0-pos, line width=5pt](3,2) .. controls (3.6,2) and (4,2.2) .. (4,3) -- (4,4);
\end{scope}
\fill[generator-2-1-1-zer] (3,2) circle (0.14);
\end{tikzpicture}}\end{aligned}
\]
\item invertible 3-cells known as \emph{cusps} which relate the unit and counit:
\[
\begin{aligned}\scalebox{0.4}{\begin{tikzpicture}
\definecolor{generator-1-0-0-pos}{RGB}{192, 57, 43}
\definecolor{generator-2-1-0-neg}{RGB}{166, 105, 8}
\definecolor{generator-2-1-0-pos}{RGB}{243, 156, 18}
\definecolor{generator-0-0-0-pos}{RGB}{41, 128, 185}
\definecolor{generator-2-1-1-zer}{RGB}{251, 221, 173}
\begin{scope}
% Background surfaces
\fill[generator-0-0-0-pos] (0,0) -- (6,0) -- (6,3) .. controls (6,3.8) and (5.6,4) .. (5,4) .. controls (4.4,4) and (4,3.8) .. (4,3) .. controls (4,2.2) and (3.6,2) .. (3,2) .. controls (2.4,2) and (2,2.2) .. (2,3) -- (2,6) -- (0,6) -- (0,0);
\fill[generator-1-0-0-pos] (6,0) -- (8,0) -- (8,6) -- (2,6) -- (2,3) .. controls (2,2.2) and (2.4,2) .. (3,2) .. controls (3.6,2) and (4,2.2) .. (4,3) .. controls (4,3.8) and (4.4,4) .. (5,4) .. controls (5.6,4) and (6,3.8) .. (6,3) -- (6,0);
% Wire layers
\draw[color=generator-2-1-0-neg, line width=5pt](3,2) .. controls (3.6,2) and (4,2.2) .. (4,3) .. controls (4,3.8) and (4.4,4) .. (5,4);
\draw[color=generator-2-1-0-pos, line width=5pt](6,0) -- (6,3) .. controls (6,3.8) and (5.6,4) .. (5,4)(3,2) .. controls (2.4,2) and (2,2.2) .. (2,3) -- (2,6);
\end{scope}
\fill[generator-2-1-1-zer] (3,2) circle (0.14);
\fill[generator-2-1-1-zer] (5,4) circle (0.14);
\end{tikzpicture}}\end{aligned}
\quad\to\quad
\begin{aligned}\scalebox{0.4}{\begin{tikzpicture}
\definecolor{generator-0-0-0-pos}{RGB}{41, 128, 185}
\definecolor{generator-1-0-0-pos}{RGB}{192, 57, 43}
\definecolor{generator-2-1-0-pos}{RGB}{243, 156, 18}
\begin{scope}
% Background surfaces
\fill[generator-0-0-0-pos] (0,0) -- (2,0) -- (2,2) -- (0,2) -- (0,0);
\fill[generator-1-0-0-pos] (2,0) -- (4,0) -- (4,2) -- (2,2) -- (2,0);
% Wire layers
\draw[color=generator-2-1-0-pos, line width=5pt](2,0) -- (2,2);
\end{scope}
\end{tikzpicture}}\end{aligned}
\]
\item invertible 4-cells known as \emph{swallowtails} which relate the different cusps:
\[
\begin{aligned}\scalebox{0.4}{\begin{tikzpicture}
\definecolor{generator-2-1-0-pos}{RGB}{243, 156, 18}
\definecolor{generator-2-1-2-pos}{RGB}{247, 188, 96}
\definecolor{generator-2-1-2-neg}{RGB}{243, 154, 13}
\definecolor{generator-0-0-1-pos}{RGB}{56, 150, 211}
\definecolor{generator-2-1-1-zer}{RGB}{251, 221, 173}

\newcommand{\wire}[2]{
  \ifdefined\recolor\draw[color=\recolor, line width=10pt]\else\draw[color=#1, line width=5pt]\fi #2
}
\newcommand{\clipped}[3]{
\begin{scope}
  \newcommand{\recolor}{#1}
  \clip#3;
  #2
\end{scope}
}

\begin{scope}[transparency group]
% Background surfaces
\fill[generator-0-0-1-pos] (0,0) -- (2,0) -- (2,3) .. controls (2,3.8) and (2.4,4) .. (3,4) .. controls (2.4,4) and (2,4.2) .. (2,5) -- (2,8) -- (0,8) -- (0,0);
\fill[generator-2-1-0-pos] (2,0) -- (8,0) -- (8,8) -- (2,8) -- (2,5) .. controls (2,4.2) and (2.4,4) .. (3,4) .. controls (2.4,4) and (2,3.8) .. (2,3) -- (2,0);
\newcommand{\layer}[1]{
  \clipped{generator-0-0-1-pos}{#1}{(0,0) -- (2,0) -- (2,3) .. controls (2,3.8) and (2.4,4) .. (3,4) .. controls (2.4,4) and (2,4.2) .. (2,5) -- (2,8) -- (0,8) -- (0,0)}
  \clipped{generator-2-1-0-pos}{#1}{(2,0) -- (8,0) -- (8,8) -- (2,8) -- (2,5) .. controls (2,4.2) and (2.4,4) .. (3,4) .. controls (2.4,4) and (2,3.8) .. (2,3) -- (2,0)}
  #1
}

% Wire layers
\wire{generator-2-1-1-zer}{(2,3) .. controls (2,3.8) and (2.4,4) .. (3,4) .. controls (3.6,4) and (4,4.2) .. (4,5) .. controls (4,5.8) and (4.4,6) .. (5,6)(6,3) -- (6,5) .. controls (6,5.8) and (5.6,6) .. (5,6)};
\layer{
\wire{generator-2-1-1-zer}{(2,0) -- (2,3)(5,2) .. controls (4.4,2) and (4,2.2) .. (4,3) .. controls (4,3.8) and (3.6,4) .. (3,4) .. controls (2.4,4) and (2,4.2) .. (2,5) -- (2,8)(5,2) .. controls (5.6,2) and (6,2.2) .. (6,3)};
}
\end{scope}
\fill[generator-2-1-2-pos] (5,2) circle (0.14);
\fill[generator-2-1-2-neg] (5,6) circle (0.14);
\end{tikzpicture}}\end{aligned}
\quad\to\quad
\begin{aligned}\scalebox{0.4}{\begin{tikzpicture}
\definecolor{generator-0-0-1-pos}{RGB}{56, 150, 211}
\definecolor{generator-2-1-0-pos}{RGB}{243, 156, 18}
\definecolor{generator-2-1-1-zer}{RGB}{251, 221, 173}
\begin{scope}
% Background surfaces
\fill[generator-0-0-1-pos] (0,0) -- (2,0) -- (2,2) -- (0,2) -- (0,0);
\fill[generator-2-1-0-pos] (2,0) -- (4,0) -- (4,2) -- (2,2) -- (2,0);
% Wire layers
\draw[color=generator-2-1-1-zer, line width=5pt](2,0) -- (2,2);
\end{scope}
\end{tikzpicture}}\end{aligned}
\]
where the source is following composite, and the target is the identity of the unit:
\[
\begin{aligned}\scalebox{0.25}{\begin{tikzpicture}
\definecolor{generator-1-0-0-pos}{RGB}{192, 57, 43}
\definecolor{generator-2-1-0-neg}{RGB}{166, 105, 8}
\definecolor{generator-2-1-0-pos}{RGB}{243, 156, 18}
\definecolor{generator-0-0-0-pos}{RGB}{41, 128, 185}
\definecolor{generator-2-1-1-zer}{RGB}{251, 221, 173}
\begin{scope}
% Background surfaces
\fill[generator-0-0-0-pos] (0,0) -- (6,0) -- (6,4) -- (4,4) -- (4,3) .. controls (4,2.2) and (3.6,2) .. (3,2) .. controls (2.4,2) and (2,2.2) .. (2,3) -- (2,4) -- (0,4) -- (0,0);
\fill[generator-1-0-0-pos] (2,3) .. controls (2,2.2) and (2.4,2) .. (3,2) .. controls (3.6,2) and (4,2.2) .. (4,3) -- (4,4) -- (2,4) -- (2,3);
% Wire layers
\draw[color=generator-2-1-0-neg, line width=5pt](3,2) .. controls (3.6,2) and (4,2.2) .. (4,3) -- (4,4);
\draw[color=generator-2-1-0-pos, line width=5pt](3,2) .. controls (2.4,2) and (2,2.2) .. (2,3) -- (2,4);
\end{scope}
\fill[generator-2-1-1-zer] (3,2) circle (0.14);
\end{tikzpicture}
}\end{aligned}
\quad\to\quad
\begin{aligned}\scalebox{0.25}{\begin{tikzpicture}
\definecolor{generator-1-0-0-pos}{RGB}{192, 57, 43}
\definecolor{generator-2-1-0-neg}{RGB}{166, 105, 8}
\definecolor{generator-2-1-0-pos}{RGB}{243, 156, 18}
\definecolor{generator-0-0-0-pos}{RGB}{41, 128, 185}
\definecolor{generator-2-1-1-zer}{RGB}{251, 221, 173}
\begin{scope}
% Background surfaces
\fill[generator-0-0-0-pos] (0,0) -- (10,0) -- (10,8) -- (8,8) -- (8,3) .. controls (8,2.2) and (7.6,2) .. (7,2) .. controls (6.4,2) and (6,2.2) .. (6,3) -- (6,5) .. controls (6,5.8) and (5.6,6) .. (5,6) .. controls (4.4,6) and (4,5.8) .. (4,5) .. controls (4,4.2) and (3.6,4) .. (3,4) .. controls (2.4,4) and (2,4.2) .. (2,5) -- (2,8) -- (0,8) -- (0,0);
\fill[generator-1-0-0-pos] (6,3) .. controls (6,2.2) and (6.4,2) .. (7,2) .. controls (7.6,2) and (8,2.2) .. (8,3) -- (8,8) -- (2,8) -- (2,5) .. controls (2,4.2) and (2.4,4) .. (3,4) .. controls (3.6,4) and (4,4.2) .. (4,5) .. controls (4,5.8) and (4.4,6) .. (5,6) .. controls (5.6,6) and (6,5.8) .. (6,5) -- (6,3);
% Wire layers
\draw[color=generator-2-1-0-neg, line width=5pt](7,2) .. controls (7.6,2) and (8,2.2) .. (8,3) -- (8,8)(3,4) .. controls (3.6,4) and (4,4.2) .. (4,5) .. controls (4,5.8) and (4.4,6) .. (5,6);
\draw[color=generator-2-1-0-pos, line width=5pt](7,2) .. controls (6.4,2) and (6,2.2) .. (6,3) -- (6,5) .. controls (6,5.8) and (5.6,6) .. (5,6)(3,4) .. controls (2.4,4) and (2,4.2) .. (2,5) -- (2,8);
\end{scope}
\fill[generator-2-1-1-zer] (7,2) circle (0.14);
\fill[generator-2-1-1-zer] (3,4) circle (0.14);
\fill[generator-2-1-1-zer] (5,6) circle (0.14);
\end{tikzpicture}
}\end{aligned}
\quad\to\quad
\begin{aligned}\scalebox{0.25}{\begin{tikzpicture}
\definecolor{generator-1-0-0-pos}{RGB}{192, 57, 43}
\definecolor{generator-2-1-0-neg}{RGB}{166, 105, 8}
\definecolor{generator-2-1-0-pos}{RGB}{243, 156, 18}
\definecolor{generator-0-0-0-pos}{RGB}{41, 128, 185}
\definecolor{generator-2-1-1-zer}{RGB}{251, 221, 173}
\begin{scope}
% Background surfaces
\fill[generator-0-0-0-pos] (0,0) -- (10,0) -- (10,8) -- (8,8) -- (8,5) .. controls (8,4.2) and (7.6,4) .. (7,4) .. controls (6.4,4) and (6,4.2) .. (6,5) .. controls (6,5.8) and (5.6,6) .. (5,6) .. controls (4.4,6) and (4,5.8) .. (4,5) -- (4,3) .. controls (4,2.2) and (3.6,2) .. (3,2) .. controls (2.4,2) and (2,2.2) .. (2,3) -- (2,8) -- (0,8) -- (0,0);
\fill[generator-1-0-0-pos] (2,3) .. controls (2,2.2) and (2.4,2) .. (3,2) .. controls (3.6,2) and (4,2.2) .. (4,3) -- (4,5) .. controls (4,5.8) and (4.4,6) .. (5,6) .. controls (5.6,6) and (6,5.8) .. (6,5) .. controls (6,4.2) and (6.4,4) .. (7,4) .. controls (7.6,4) and (8,4.2) .. (8,5) -- (8,8) -- (2,8) -- (2,3);
% Wire layers
\draw[color=generator-2-1-0-neg, line width=5pt](3,2) .. controls (3.6,2) and (4,2.2) .. (4,3) -- (4,5) .. controls (4,5.8) and (4.4,6) .. (5,6)(7,4) .. controls (7.6,4) and (8,4.2) .. (8,5) -- (8,8);
\draw[color=generator-2-1-0-pos, line width=5pt](3,2) .. controls (2.4,2) and (2,2.2) .. (2,3) -- (2,8)(7,4) .. controls (6.4,4) and (6,4.2) .. (6,5) .. controls (6,5.8) and (5.6,6) .. (5,6);
\end{scope}
\fill[generator-2-1-1-zer] (3,2) circle (0.14);
\fill[generator-2-1-1-zer] (7,4) circle (0.14);
\fill[generator-2-1-1-zer] (5,6) circle (0.14);
\end{tikzpicture}
}\end{aligned}
\quad\to\quad
\begin{aligned}\scalebox{0.25}{\begin{tikzpicture}
\definecolor{generator-1-0-0-pos}{RGB}{192, 57, 43}
\definecolor{generator-2-1-0-neg}{RGB}{166, 105, 8}
\definecolor{generator-2-1-0-pos}{RGB}{243, 156, 18}
\definecolor{generator-0-0-0-pos}{RGB}{41, 128, 185}
\definecolor{generator-2-1-1-zer}{RGB}{251, 221, 173}
\begin{scope}
% Background surfaces
\fill[generator-0-0-0-pos] (0,0) -- (6,0) -- (6,4) -- (4,4) -- (4,3) .. controls (4,2.2) and (3.6,2) .. (3,2) .. controls (2.4,2) and (2,2.2) .. (2,3) -- (2,4) -- (0,4) -- (0,0);
\fill[generator-1-0-0-pos] (2,3) .. controls (2,2.2) and (2.4,2) .. (3,2) .. controls (3.6,2) and (4,2.2) .. (4,3) -- (4,4) -- (2,4) -- (2,3);
% Wire layers
\draw[color=generator-2-1-0-neg, line width=5pt](3,2) .. controls (3.6,2) and (4,2.2) .. (4,3) -- (4,4);
\draw[color=generator-2-1-0-pos, line width=5pt](3,2) .. controls (2.4,2) and (2,2.2) .. (2,3) -- (2,4);
\end{scope}
\fill[generator-2-1-1-zer] (3,2) circle (0.14);
\end{tikzpicture}
}\end{aligned}
\]
\end{itemize}
\Cref{fig:cusp-3d} shows the cusp rendered in 3D, and \Cref{fig:swallowtail-source-3d} shows the source of the swallowtail rendered in 3D.
The swallowtail itself can be rendered in 4D; this can be viewed on YouTube:
\begin{center}
\url{https://youtu.be/UZ6mBlbrpPc}
\end{center}

\footnotetext{There is an infinite sequence of these coherences, which are the study of catastrophe theory \cite{Zee76}.}

\end{document}